\documentclass{amsart}
\usepackage[T1]{fontenc}
\usepackage{amssymb}
\usepackage[hyperref]{degt}

\def\2{\color{red}}
\def\3{\color{green}}
\def\4{\color{blue}}

\def\gf#1#2{[#1_{#2}]}

\def\KK-{\pdfstr{K3}{$K3$}-\penalty0\ignorespaces}

\let\0\varnothing

\def\dual{^\vee}

\let\graph\Gamma
\let\pencil\Pi
\let\fiber\Sigma
\let\Weyl\Delta

\def\Fano{\Cal{F}}

\def\Fn{\operatorname{Fn}}
\def\Fnex{\Fn^\mathrm{ex}}
\def\extended{^{\mbox{\tiny ex}}}
\def\extended{^\mathrm{ex}}

\def\graphex{\graph\extended}

\def\root{\operatorname{root}}
\def\sat{\operatorname{sat}}

\def\val{\operatorname{val}}
\def\Star{\operatorname{star}}

\def\NS{\operatorname{NS}}

\def\Nef{\operatorname{Nef}}

\def\kappal{\mathfrak{l}}

\def\graphex{\graph\extended}

\def\resp.{resp\PERIOD}

\let\M=M
\let\o=o

\let\ex=e
\let\iso=p

\def\CK{\Cal K}

\def\bL{\bold{L}}

\def\CC{\Cal C}

\def\bl{\bar l}

\def\fm{\frak{m}}
\def\fd{\frak{d}}
\def\order{\operatorname{order}}

\def\Sym{\operatorname{Sym}}

\def\tto{\mathrel{\to}}
\def\iin{\mathrel{\subset}}

{
\catcode`\1\active
\catcode`\.\active
\gdef\tconfig{\vcenter\bgroup\offinterlineskip
\def\1{&\bullet}\def\.{&\cdot}\def\2{&\circ}%
\halign\bgroup
 \hbox to0pt{\hss##\hss}&&\hbox to9pt{\small\hss$\mathstrut##$\hss}\cr}
\gdef\endtconfig{\crcr\egroup\egroup}
}

\def\biquad{\Omega}

\def\itemrefform#1{\ifmmode(#1)\else$(#1)$\fi}

\def\lineno{\afterassignment\linenoii\count0=}
\def\linenoii{\let\next=L\ifnum\count0>12\let\next=M\advance\count0-12\fi
 \next_{\the\count0}}
\def\DEFS{\let\*\times\def\"{^*}\let\l\lineno\let\c\mathbf}%

\def\i{^\mathrm{i}}
\def\ii{^\mathrm{ii}}

\def\HLform#1{\{#1\}}
\let\HLcomma,
\def\HL[#1]{\HLform{\let\comma\relax\HLL#1,\end,}}
\def\HLL#1,{\ifx#1\end\let\comma\relax\else\ln#1\expandafter\HLL\fi}
\def\ln{\afterassignment\lnn\count0=}
\def\lnn{\comma\ifnum\count0>12\advance\count0 -12\relax m\else l\fi
 _{\the\count0}\let\comma\HLcomma}

\def\Humbert{\mathrm{H}}
\def\bd{\bar{d}}
\def\bl{\bar{l}}

\def\same{\simeq}

\addtheorem{Algorithm}{algorithm}
\newcounter{graph}
\def\thegraph{\arabic{graph}}

\title{Lines on \KK-sextics with simple singularities}

\author{Alex Degtyarev}

\address{%
Department of Mathematics\\
Bilkent University\\
06800 Ankara, TURKEY}

\email{
degt@fen.bilkent.edu.tr}

\author{S\l awomir Rams}
\address{Theoretical Computer Science Department, Faculty of Mathematics and Computer Science, Jagiellonian University,
	ul.\ {\L}ojasiewicza 6,  30-348 Krak\'{o}w, Poland}
\email{slawomir.rams@uj.edu.pl}

\thanks{%
{A.D.} was partially supported by the T\"{U}B\DOTaccent{I}TAK grant 123F111%
}

\thanks{%
	S.R. was supported by the National Science Centre, Poland, Opus  grant
	no.\ 2024/\allowbreak53/\allowbreak B/\allowbreak ST1/\allowbreak00161 and
	partially supported by
	{T\"{U}B\DOTaccent{I}TAK B\DOTaccent{I}DEB
		2221 Visiting Scientist Fellowship Program}%
}



\begin{document}

\begin{abstract}
We advance our understanding of the configurations of low degree smooth
rational curves on (quasi-)polarized complex $K3$-surfaces.
We apply
our efficient approach
to classify the configurations of at least $36$ lines
on $K3$-sextics with at
worst
A--D--E singularities.
Besides, we
characterize a certain class of infinite dihedral groups of birational
automorphisms of $K3$-sextics,
we show that no $K3$-sextic can contain a Kummer configuration of lines,
and we give a
complete
account of the line configurations on
the closest
analogue of Kummer $K3$-octics
or quartics, \viz. the so-called Humbert $K3$-sextics.
\end{abstract}

\maketitle

\section{Introduction}\label{S.intro}

We work over the field $\C$ of complex numbers.

\subsection{Motivation}\label{s.motivation}
This work is inspired by the recent paper~\cite{DDK} and by
our previous study of large configurations of lines on $K3$-complete
intersections with $\bA$--$\bD$--$\bE$ singularities, \viz.
quartics~\cite{Degtyarev.Rams.triangular} and octics~\cite{degt.Rams:octics}.
In~\cite{DDK}, a fairly complete description of the unexpectedly rich geometry
of $K3$-sextics with an abstract configuration $(12_6,12_6)$ of lines,
the so-called \emph{Humbert sextics}, is presented; here, we show that,
although quite far from the maximum (which is $42$ lines, see \autoref{th.main}), this
configuration of $24$ lines is remarkable in many ways.

To begin with, we observe that, unlike quartics or octics, a $K3$-sextic cannot contain a  Kummer configuration of lines; in
fact, the maximal number of pairwise skew lines on a sextic is $12$, see
\autoref{th.Kummer}. Humbert is the only configuration of lines on a smooth
sextic that has two such $12$-tuples.

Next, recall (see \cite{DDK} or \autoref{s.tau} below)
that a line $\ell\subset X$ on a $K3$-sextic gives rise to a two-to-one map
$X\to\Cp2$ and, hence, to the
\emph{deck translation} $\tau\:X\to X$, which is a birational automorphism
of~$X$. On a Humbert sextic, the deck translations defined by any two
distinct lines do not commute; moreover, they generate an infinite dihedral
group. We extend this statement to \emph{almost} any pair of lines on a
$K3$-sextic, see \autoref{th.aut} below and references therein. The ``almost'' here
is crucial: the action of $\tau$ depends on the other lines and conics
intersecting~$\ell$, so that \autoref{th.aut} is not a mere application of
simple linear algebra. The (relatively) few exceptions are pairs of lines
appearing in quite sophisticated explicit configurations.

The space of Humbert sextics has strata with more than the original $24$
lines; this stratification was announced in \cite[Remark 7.1]{DDK} and used
in some statements in \latin{loc.~cit}, \eg, in the proof of Lemma~5.1.
We explain the details and complete
the proof in the present paper, see \autoref{S.Humbert}. 

Finally, recall that, according to Saint-Donat, the space of $K3$-octics in
$\Cp5$ has a codimension~$1$ stratum of \emph{special} octics, \ie, those
that are not triquadrics; homologically, this stratum is characterized by the
presence of so-called \emph{$3$-isotropic vectors} in the N\'eron--Severi
lattice. We discover that, in the realm of sextics, the presence of such a
vector is also geometrically meaningful (see \autoref{def.special}); as in
the case of octics, it affects the configurations of lines, \cf.,
\eg, \autoref{cor.boundsforspecial}.
The Humbert configuration is the maximal configuration of lines that may
appear on both special and non-special sextics.

According to these goals, the paper splits naturally into two parts. First,
we develop techniques that enable us to classify $K3$-sextics with large
configurations of lines, culminating in a proof of \autoref{th.main}.
Then, we apply these techniques to the study of more special, not necessarily
very large, configurations, including a neighborhood of a fixed pair of lines and
the Humbert configuration $(12_6,12_6)$.

%
%

\subsection{Main results}



Recall that the configuration
of lines on a surface $X$ is encoded in its  \emph{Fano graph}
$\Fn X$ (see \autoref{s.LLL} for details) and
an \emph{equilinear stratum}  of quasi-polarized \KK-sextics is defined as an  irreducible component of the space
\[*
\bigl\{\text{$X\tto\Cp4$ a $K3$-sextic}\bigm|
\ls|\Fn X|=\const,\
\mu(X)=\const\bigr\} \, ,
\]
where $\mu$ is the total Milnor number of the singularities
of~$X$.
We
prove the following theorem
(see \autoref{Ss:Fanographs} below for the precise description of the
contents of \autoref{tab.list}).
We recall that a \KK-surface $X$ {in} $\Cp{n}$ may have at worst
$\bA$--$\bD$--$\bE$ singularities.

\theorem[see \autoref{proof.main}]\label{th.main}
Let $X \tto \Cp4$ be a
sextic \KK-surface.
Then, $X$ contains at most $42$ lines {\rm(}at most $36$ if
$\Sing X \neq \emptyset${\rm)}. Futhermore, if $X$ contains at least
$36$ lines, then it belongs to one of
the
$30$
equilinear strata listed in   \autoref{tab.list}.
\table
\caption{$K3$-sextics with at least $36$ lines
 (see \autoref{th.main} and \autoref{s.notation})}\label{tab.list}
\hbox to\hsize{\hss\vbox{\offinterlineskip%
		\def\config#1(#2){#1}%
		\def\ps#1{}%
		\def\aut#1{#1}%
		\def\symplectic(#1,#2){\def\symsize{#1}\def\symindex{#2}}%
		\def\group#1(#2){\symsize\rlap{$_{\symindex}^{\ifnum#2>1 #2\fi}$}}%
		\def\group#1(#2){(\symsize\mathinner:\symindex)\ifnum#2>1 \rlap{$^{#2}$}\fi}%
		\def\mat#1{#1}%
		\def\mmat#1{#1}%
		\def\*{\llap{$^*$}}%
		\def\bU(#1){\bold{U}_{#1}}%
		\def\bU(#1){\bold{U}(#1)}%
		\let\+\oplus
		\def\rc#1{#1}%
		\def\i{'}\def\ii{''}\def\iii{'''}\def\0{}\def\v{^{\mathrm{v}}}%
\def\([#1,#2]){\refstepcounter{graph}\gf{\thegraph}{#1}}
		\def\3{\rlap{$^3$}}%
		\halign{\strut\quad\hss$#$\hss\quad&#&#&\hss$#$\hss\quad&\hss$#$\hss\quad&\hss$#$\hss\quad&
			\ \hss$#$\hss\quad&\ $#$\hss\quad\cr
			\noalign{\hrule\vskip2pt}
			\graphex&&&\Sing X_6&\ls|\Aut\graph|&\ls|\Sym X_6|&(r,c)&\NS(X_6)^\perp\cr
			\noalign{\vskip1pt\hrule\vskip2pt}
			\config \Psi_{42}(19)
			\ps{(6,3)^{2}}&&&
			&\aut{432} &\symplectic(36,9) \group72(2)
			&\rc{(1,0)} &\mmat{\*{-\bA_2}\+[-18]}\cr
			\config \Psi_{38}(19)
			\ps{(5,4)^{2}}&&&
			&\aut{32} &\symplectic(4,1) \group8(2)
			&\rc{(1,0)} &\mmat{\*{-\bA_2}\+[-24]}\cr
			\config \Psi_{36}'(18)
			\ps{(6,0)^{2}}&&&
			&\aut{1296} &\symplectic(18,4) \group36(2)
			&\rc{(1,0)} &\mmat{\*{-\bA_2}\+\bA_2(3)}\cr
			\config \Psi_{36}''(19)
			\ps{(6,3)^{1} (3,6)^{1}}&&&
			&\aut{144} &\symplectic(6,1) \group12(2)
			&\rc{(1,0)} &\mmat{\*{-\bA_2}\+[-36]}\cr
			\config \([ 36, 0 ])(19)
			\ps{}&&&
			&\aut{144} &\symplectic(18,4) \group36(2)
			&\rc{(1,0)} &\mmat{\*{-\bA_2}\+[-24]}\cr
			\config \Theta^\circ_{36}(20)
			\ps{}&&&
			&\aut{48} &\symplectic(24,12) \group24(1)
			&\rc{(1,0)} &\mat{\*[2,0,36]}\cr
			\config \([ 36, 0 ])(20)
			\ps{}&&&
			&\aut{48} &\symplectic(12,3) \group24(2)
			&\rc{(1,0)} &\mat{[4,2,14]}\cr
			\config \Theta_{36}(19)
			\ps{}&&&
			&\aut{24} &\symplectic(6,1) \group12(2)
			&\rc{(1,0)} &\mmat{\*[16]\+[2,1,-2]}\cr
			\config \([ 36, 0 ])(20)
			\ps{}&&&
			&\aut{16} &\symplectic(4,2) \group8(2)
			&\rc{(0,1)} &\mat{[6,2,10]}\cr
			\config \([ 36, 0 ])(20)
			\ps{}&&&
			&\aut{8} &\symplectic(2,1) \group4(2)
			&\rc{(1,0)} &\mat{[6,3,14]}\cr
			\config \([ 36, 0 ])(20)
			\ps{}&&&
			&\aut{4} &\symplectic(2,1) \group4(2)
			&\rc{(1,0)} &\mat{\*[2,1,30]}\cr&&&&&
			&\rc{(0,1)} &\mat{[6,1,10]}\cr
			\config \([ 36, 3 ])(20)
			\ps{}&&&3\bA_1
			&\aut{48} &\symplectic(12,4) \group24(2)
			&\rc{(1,0)} &\mat{\*[2,0,16]}\cr
			\config \([ 36, 3 ])(20)
			\ps{}&&&3\bA_1
			&\aut{36} &\symplectic(18,4) \group36(2)
			&\rc{(1,0)} &\mat{\*[2,1,14]}\cr
			\config \([ 36, 2 ])(20)
			\ps{}&&&2\bA_1
			&\aut{96} &\symplectic(24,12) \group48(2)
			&\rc{(1,0)} &\mat{[6,2,6]}\cr
			\config \([ 36, 2 ])(20)
			\ps{}&&&2\bA_1
			&\aut{24} &\symplectic(6,1) \group12(2)
			&\rc{(1,0)} &\mat{[6,1,6]}\cr
			\config \([ 36, 2 ])(20)
			\ps{}&&&2\bA_1
			&\aut{16} &\symplectic(4,2) \group8(2)
			&\rc{(1,0)} &\mat{[4,0,10]}\cr
			\config \([ 36, 2 ])(20)
			\ps{}&&&2\bA_1
			&\aut{16} &\symplectic(4,2) \group8(2)
			&\rc{(1,0)} &\mat{\*[2,0,18]}\cr
			\config \Psi_{36}^1(20)
			\ps{(5,3)^{2}}&&&\bA_1
			&\aut{8} &\symplectic(2,1) \group4(2)
			&\rc{(1,0)} &\mat{\*[2,0,24]}\cr
			\config \Theta^1_{36}(20)
			\ps{}&&&\bA_1
			&\aut{24} &\symplectic(6,1) \group12(2)
			&\rc{(2,0)} &\mat{[10,0,16]}\cr
			\config \([ 36, 1 ])(20)
			\ps{}&&&\bA_1
			&\aut{12} &\symplectic(6,1) \group12(2)
			&\rc{(1,0)} &\mat{\*[2,1,22]}\cr
			\config \([ 36, 1 ])(20)
			\ps{}&&&\bA_1
			&\aut{8} &\symplectic(4,2) \group8(2)
			&\rc{(1,0)} &\mat{\*\*[2,0,22]}\cr&&&&&
			&\rc{(0,1)} &\mat{[6,2,8]}\cr
			\config \([ 36, 1 ])(20)
			\ps{}&&&\bA_1
			&\aut{8} &\symplectic(2,1) \group4(2)
			&\rc{(1,0)} &\mat{\*[2,1,26]}\cr
			\config \([ 36, 1 ])(20)
			\ps{}&&&\bA_1
			&\aut{8} &\symplectic(2,1) \group4(2)
			&\rc{(1,0)} &\mat{\*[2,0,24]}\cr
			\noalign{\vskip1pt\hrule}%
			\crcr}}\hss}%
\endtable
\endtheorem

Let us emphasize that the powerful tool developed in \cite{degt:lines,degt.Rams:octics}
is strong enough to classify line configurations
on quasi-polarized \KK-surfaces of degree
eight or higher, and we do not expect anything new (as compared
to~\cite{degt:lines}), except, probably, a confirmation of the conjecture
that the presence of exceptional divisors reduces the maximal number of
lines.
However,
extra ideas are required to deal with triangular configurations
on quartics (see \cite{Degtyarev.Rams.triangular}) and certain $\tA_3$-configurations on
sextics.
It is this latter case,
\ie, \emph{quadrangular} Fano graphs,
that constitutes the bulk of the proof of \autoref{th.main}
(see \autoref{s.biquads}) and all but two configurations in
\autoref{tab.list}.

Our proof of \autoref{th.main} is considerably different from that of its
analogue for octics, see~\cite{degt.Rams:octics}. Whereas in the former case
a substantial technical difficulty (\viz. large geometric kernels,
see \cite[\S\,6.2]{degt.Rams:octics})
was presented by Kummer and almost Kummer
surfaces, here we show (see \autoref{th.Kummer}) that a \KK-sextic may have
at most~$12$ pairwise disjoint lines, so that Kummer sextics are of limited
relevance. Instead, we concentrate on the quadrangular graphs,
see~\autoref{S.quad}. (As we explained above, it appears that the closest analogues of Kummer
surfaces are the \emph{Humbert sextics}~\cite{DDK}; they are discussed in
\autoref{S.Humbert}, where we also illustrate our approach by more explicit
examples.)

As a similarity between the two cases, we discover a natural subdivision
of \KK-sextics into \emph{special} and \emph{non-special}, see
\autoref{def.special}. As in the case of octics, this property is detected
by the presence of a $3$-isotropic vector in the
N\'{e}ron--Severi lattice (see \autoref{th.special}) and,
as in \cite[Theorem~1.3 and Proposition~5.9]{degt:lines} (see also
\cite[Theorem~1.2]{degt.Rams:octics}), it affects dramatically the
configurations of lines.

\corollary[see \autoref{s.special}]\label{cor.boundsforspecial}
The maximal number of lines on a non-special {\rm(}resp.\ special{\rm)}
\KK-sextic
$X$
is $42$ {\rm(}resp.\ $36${\rm)}.
If $\Sing X \neq \emptyset$, the maximum number of lines  is $36$
in both cases.
A
maximal configuration realized by both special and non-special sextics
has $24$ lines, \eg, that of Humbert sextics
{\rom(see \autoref{rem.Humbert.spec}\rom).}
\endcorollary

As stated in \autoref{s.motivation}, we go far
beyond a mere classification of large configurations
and illustrate the power
of our approach
and its applicability
to a number of other problems, not necessarily related to line counting.
In particular, we settle a few questions raised in~\cite{DDK}.
Most
notably, we find a plethora of \emph{explicit} infinite subgroups of the
groups $\Aut X$ of abstract (\latin{a.k.a}.\ birational) automorphisms of \KK-surfaces.
(For the definition of the deck translation $\tau\:X\to X$
corresponding to a line~$\ell$ on a \KK-sextic~$X$, see \autoref{s.tau}).

\theorem[see \autoref{proof.aut}]\label{th.aut}
Let $X\iin\Cp4$ be a smooth $K3$-sextic, $\ell\i\ne\ell\ii$ a pair of
lines on~$X$, and $\tau\i,\tau\ii\:X\to X$ the corresponding deck
translations. Then, with the \rom(very few\rom) exceptions classified in
Lemmata~\ref{lem.aut.intr},~\ref{lem.aut.skew}\rom,
the involutions
$\tau\i$ and $\tau\ii$ generate
an infinite dihedral group.
\endtheorem

\corollary[see \autoref{proof.cor.aut}]\label{cor.aut}
If a smooth $K3$-sextic $X\iin\Cp4$ contains
at least two lines, then the group $\Aut X$ of birational automorphisms
of~$X$
is infinite. Moreover, there is a pair of lines on~$X$ whose deck
translations generate an infinite dihedral group.

In particular, if a \KK-sextic contains at least two lines, then it
contains infinitely many \emph{smooth} rational curves {\rm(}\cf.~\autoref{ex.aut}{\rm)}.
\endcorollary

Apart from the existence of smooth rational curves,
our interest in compositions of deck translations stems from the fact that
\autoref{cor.aut} provides
over 9000 fairly explicit examples of automorphisms
of projective \KK-surfaces with non-trivial entropy.
The last assertion of \autoref{cor.aut} is illustrated by the following
example, interesting since very little is known about high degree smooth
rational curves on $K3$-quartics, sextics, or octics.

\example\label{ex.aut}
The smooth rational curves given by \autoref{cor.aut} may be
hard to find! Using
Vinberg's algorithm (see \autoref{s.LLL}), we found
that, if
the  group
$\NS(X)$ is generated by the class $h$ of the degree-$6$ polarization and  two skew lines $\ell\i\ne\ell\ii$,
then $X$ contains but $13$
smooth rational curves of degree up to $3000$: the lines $\ell\i,\ell\ii$,
the quartic $h-\ell\i-\ell\ii$, and a ``symmetric'' pair of curves
in each degree
\[*
14, 41, 164, 601, 2254.
\]
In a sense, this problem is related to
Pell type equations, whose solutions may be very scarce.
It appears that, in this case, $\Aut X=\DG\infty$ is generated by
$\tau\i$, $\tau\ii$
and that the smooth rational curves
constitute a single $\DG\infty$-orbit.

If $\ell\i,\ell\ii$ above intersect, there are more curves, two in each
degree
\[*
1, 13, 37, 73, 121, 181, 253, 337, 433, 541, 661, 793, 937,\ldots
\]
They seem to constitute two $\DG\infty$-orbits.
\endexample

Another application is a complete description of the
equilinear stratification of the space of Humbert sextics. Precise statements
are postponed till \autoref{S.Humbert}, as they need a great deal of
preparation. One consequence of this description is the fact that the
``original'' configuration of $12+12$ lines is preserved by the
\emph{projective} automorphisms of any Humbert sextic.

\subsection{Notation in \autoref{tab.list}} \label{Ss:Fanographs}
As stated in \autoref{th.main}, the complete classification of configurations of at least $36$ lines on quasi-polarized \KK-surfaces is presented in \autoref{tab.list},
where we collect the following data:
\roster*
\item
the name of the extended Fano graph~$\graphex$,
the subscript always referring to the number of lines;
$\Psi_*^*$ are triangular graphs, others are quadrangular;
\item
the size of the group $\Aut\graphex$ of abstract graph automorphisms of~$\graphex$;
\item
the group $\Sym X$ of symplectic automorphisms of a generic
sextic~$X$ with
the given Fano graph,
in the form $(\text{size}:\text{index})$, referring to the
\texttt{SmallGroup}
 library in \GAP~\cite{GAP4.13};
the superscript is the index of $\Sym X$ in the full group $\Aut(X,h)$
of
projective automorphisms of~$X$ (if greater than~$1$);
\item
the numbers $(r,c)$ of, respectively, real and pairs of complex conjugate
equilinear strata;
\item
the (generic, if $\rank\ge3$) transcendental lattice $T:=\NS(X)^\perp$;
it is marked with a $^*$ if the corresponding
stratum has
a real quartic with all lines real (see \cite[Lemma 3.8]{DIS}).
The lattice notation is explained in \autoref{s.notation} below.
\endroster
Whenever the
isomorphism class of $T$ is not uniquely determined by the graph~$\graphex$,
each
lattice is listed in a separate row (following the main entry),
and the numbers $(r,c)$ of components are itemized accordingly.

\remark \label{rem:discussion.graphs}
(1)
Formally, in the first column of \autoref{tab.list} we present
\emph{extended} Fano graphs
(see \autoref{alg.Vinberg}).
However, we assert that,
with one exception, \viz. $\Theta_{32}$ and $\Theta_{32}^1$, see
\autoref{rem:lem.111}, the
corresponding
\emph{Fano graphs}
are pairwise distinct.

\smallskip
(2) Surprisingly, none   of the two extremal octic
configurations from~\cite[Table~1]{degt.Rams:octics}
appear here. Indeed,
the graph $\Theta'_{36}$ has Kummer index $14$
and, thus, it is ruled out by \autoref{th.Kummer}.
The other graph
$\Theta''_{36}$ has Kummer index $12$;
it is is admissible but fails to be
geometric for $h^2=6$.

\smallskip
(3) There are too many large quadrangular extended Fano graphs of \KK-sextics to
denote them with distinct
letters;
therefore, we merely number them.
Still, the graph $\Theta_{36}$ is of special interest  because it
appears as the solution to a few other line counting problems:
\roster*
\item
it is the maximal degeneration of Humbert sextics
(see \autoref{th.Humbert});
\item
it maximizes the number of $h$-fragments on a non-special
sextic (see \cite{degt:hplanes});
\item
it is realized by both smooth and singular $K3$-sextics, see $\Theta^1_{36}$.
\endroster
Another special graph, $\Theta^\circ_{36}$, is the only
quadrangular graph
with $36$ vertices
that has no
biquadrangles, see \autoref{lem:tA3biquadranglefree}.

\smallskip
(4) 	By  direct computation, none of the quadrangular graphs in \autoref{tab.list}
is special (see \autoref{rem.vocabulary}): smooth are not
bipartite, see \autoref{lem.bipartite}; others are of rank~$20$ and
their lattices
do not have a $3$-isotropic vector (see \autoref{th.special}).
\endremark


\subsection{Contents of the paper} \label{Ss:contents}
{\leavevmode In \autoref{S.prelim},}
we recall the necessary notation and collect various well-known facts that are
fundamental for our
work.
Then, in \autoref{s.Kummer},  we study
the so-called \emph{Kummer index} of \KK-sextics
and show that it never exceeds $12$ (see \autoref{th.Kummer}). In \autoref{S.bounds}
we prove
a few
basic bounds for line configurations on \KK-sextics and present the taxonomy of
hyperbolic graphs and pencils (\autoref{s,pencils}). After those
preparations,
 in \autoref{S.quad} we study quadrangular graphs in more detail
 and provide the crucial part of the proof of \autoref{th.main}.
 Finally, we examine triangular graphs (\autoref{lem:triangular.pencil})
and those
of other types
(merely listing the results obtained by the
standard techniques of \cite{degt.Rams:octics})
and
conclude the proof of \autoref{th.main}
in \autoref{proof.main}.

In \autoref{S.aut} we discuss
certain automorphisms of \KK-sextics,
\viz. the
deck translations corresponding to
lines, see \autoref{s.tau}.
In \autoref{s.special}
we characterize
the so-called
special sextics.
Then, upon a detailed analysis of the joint star of a pair of lines in
\autoref{s.aut.trig}--\autoref{s.aut.skew},
we prove \autoref{th.aut} and  \autoref{cor.aut} in
\autoref{proof.cor.aut}.

Finally, in \autoref{S.Humbert}, we discuss the geometry of the so-called
Humbert sextics  (see \cite{Humbert} and \cite{DDK}). This class of
\KK-sextics provides a number of interesting examples
(\eg, \cite[Corollary~4.7]{Humbert} suggests that it can be seen as a counterpart of  Kummer surfaces in the realm of sextics with many lines)
and gives us an opportunity to illustrate our algorithms by
working out the details of the computation.
In particular,
\autoref{th.Humbert} gives rise to
a purely geometric construction of one of the equilinear strata
presented in \autoref{tab.list},
\viz. $\Theta_{36}$ discussed in
\autoref{rem:discussion.graphs}.

 \subsection{Acknowledgements}
We started our  study of lines on \KK-surfaces
	 during Alex Degtyarev's
	visit to the Jagiellonian
	University (Cracow, Poland). A.D.\ is profoundly grateful to this institution for the
	hospitality and excellent working conditions.
Substantial part of
{the} paper was written during S{\l}awomir Rams' stay at
the Department of Mathematics of Bilkent University (Ankara, Turkey) within
the scope of
 T\"{U}B\DOTaccent{I}TAK B\DOTaccent{I}DEB
(2221 Visiting Scientist Fellowship Program).
S.R.\ would like to thank Bilkent University and T\"{U}B\DOTaccent ITAK for creating perfect
working conditions and numerous inspiring discussions with members of the
Department of Mathematics of Bilkent University.

\section{Preliminaries}\label{S.prelim}

In this section we collect/recall various facts
needed in the sequel.
We refer to \cite{degt.Rams:octics} for a more detailed exposition.

\subsection{Common notation}\label{s.notation}

We  maintain the standard lattice-theoretic notation (the determinant,
discriminant group,
\etc.), as introduced in \cite{Conway.Sloane,Nikulin:forms}.
Futhermore,
we follow the notation introduced in
 \cite[\S\,1.4]{degt.Rams:octics} (see also \cite{degt:lines,Degtyarev.Rams.triangular}):
\roster*
\item
$\bA_p$, $p\ge1$, $\bD_q$, $q\ge4$, $\bE_6$, $\bE_7$, $\bE_8$
are the \emph{negative definite} root
lattices generated by the corresponding indecomposable root systems
(as in~\cite{Bourbaki:Lie});
\item
$[a]:=\Z u$ is the  rank~$1$ lattice such that  $u^2=a$;
\item
$[a,b,c]:=\Z u+\Z v$, $u^2=a$, $u\cdot v=b$, $v^2=c$, is a lattice of
rank~$2$; when it is positive definite, we assume that $0<a\le c$ and
$0\le2b\le a$: then, $u$ is a shortest vector, $v$ is a next shortest one,
and the triple $(a,b,c)$ is unique;
\item
$\bU:=[0,1,0]$ is the \emph{hyperbolic plane};
\item $L(n)$ denotes the lattice obtained by the scaling of
a given lattice~$L$ by a fixed integer  $n\in\Z$,
as opposed to $nL:=L^{\oplus n}$;
however, $-L:=L(-1)$;
\item
$L\dual:=\Hom(L,\Z)$ (resp.\ $\discr L$) denotes the dual group
(resp.\ the
\emph{discriminant group} with the induced
$\Q/2\Z$-valued quadratic form~$q$; \cf.
\cite{Nikulin:forms}); 

\item
 $\sigma_{\pm,0}(L)$ stand for
the  inertia indices of the form
$L \otimes \R$.
\endroster
As in \cite{degt.Rams:octics},
all lattices considered are even (unless stated otherwise).
\emph{Isometries} are not assumed injective or surjective: they merely preserve
the quadratic form.

\subsection{Polarized lattices}\label{s.lattice}
A \emph{$2n$-polarized lattice} is a non-degenerate even hyperbolic (\ie, such
that
$\Gs_+N=1$) lattice~$N$ equipped with a distinguished vector $h\in N$ of
the prescribed square $h^2=2n>0$. The motivating example is the
N\'{e}ron--Severi lattice $\NS(X)$ of a quasi-polarized $K3$-surface
$X\to\Cp{n+1}$, in which $h$ is the class of a hyperplane section, $h^2=2n$.
Another common example is the lattice
\[
\Fano_{2n}(\graph):=(\Z\graph+\Z h)/\!\ker,\quad h^2=2n,
\label{eq.Fano}
\]
defined by a colored loop free (multi-)graph~$\graph$,
\emph{assuming that $\Gs_+\Fano_{2n}(\graph)=1$}.
Here, denoting the coloring function by $\deg\:\graph\to\{0,\ldots,d\}$, we
define
\roster*
\item
$v^2=-2$ and $v\cdot h=\deg(v)$ for a vertex $v\in\graph$ and
\item
$u\cdot v$ is the number of edges $[u,v]$ for $u,v\in\graph$.
\endroster
In this paper the default color is~$1$ (we mostly deal with lines),
the default value is $2n=6$ (sextic $K3$-surfaces), and graphs are usually
simple.

\observation\label{obs.ker}
If $N$ is a nondegenerate hyperbolic lattice, then any
sublattice $S\subset N$ with $\Gs_+S=1$ is also nondegenerate.
Hence, if $N\ni h$ is polarized,
for any isometry $\iota\:S\to N$
such that $\iota(S)\ni h$, the lattice $S$ is hyperbolic and
$\ker S=\Ker\iota$.
We often use the following simplified version:
for any $u,v\in N$, one has
\[*
\det(\Z h+\Z u+\Z v)\ge0
\]
and, if $\det=0$, the vectors $h,u,v$ are linearly dependent.
\endobservation

\algorithm\label{alg.LLL}
Many algorithms
require finding all/some vectors $v\in N$ of a given
square and degree:
\[*
v^2=s,\qquad v\cdot h=d.
\]
To this end,
we consider the negative
definite lattice $h^\perp\subset N$ and list all (finitely many) vectors
$\bar{v}\in h^\perp$ of square $\bar{v}^2=4n^2s-2nd^2$. Then, we select those
for which the \latin{a priori} rational vector $v=(dh-\bar{v})/2n$ is in~$N$.
(We do not discuss various optimization tricks here.)
The listing of all vectors $\bar{v}\in h^\perp$ is the only part of our argument
that, considering the ranks involved, is difficult to be explained/done by
hand. We use the LLL-algorithm
(Lenstra--Lenstra--Lov\'{a}sz lattice basis reduction
algorithm)
implemented as \texttt{ShortestVectors} in
\texttt{GAP}~\cite{GAP4.13}.
\endalgorithm

\subsection{Rational curves on (quasi-)polarized $K3$-surfaces}\label{s.LLL}
A
\emph{$2n$-quasi-polarized \KK-surface} is a smooth minimal \KK-surface~$X$
mapped to $\Cp{n+1}$ by  a complete base-point free linear system $|h|$, where
$h\in\NS(X)$ is a big and nef degree $2n$ class, called the \emph{quasi-polarization}. (Usually we assume that the map $X \rightarrow \Cp{n+1}$ is birational onto its image, as opposed to hyperelliptic, \cf.~\cite{Saint-Donat}.) By a \emph{line} (or,
more generally, \emph{degree~$d$ smooth rational curve}) on~$X$ we mean a
smooth rational curve $c\subset X$ whose image in~$\Cp{n+1}$ has degree~$d$.
Degree~$0$ curves are contracted to points, and the image of a degree $d\ge3$ curve
does not need to be smooth. We never work directly with the image
$X_{2n}\subset\Cp{n+1}$ of~$X$ itself. The \emph{intersection} of two
(smooth) rational curves is their intersection in~$X$.

Let $X$ be a quasi-polarized $K3$-surface and $N=\NS(X)\ni h$.
Then, the degree $d>0$
genus~$0$ curves on~$X$ are in a bijection with the set
\[*
\root_d(N,h):=\bigl\{v\in N\bigm|v^2=-2,\ v\cdot h=d\bigr\},
\]
which is easily computed using \autoref{alg.LLL}.
The case $d=0$ is special:
there is a partition
\[
\root_0(N,h)=P\cup(-P)
\label{eq.partition}
\]
into \emph{positive} and
\emph{negative} roots (see~\cite{Bourbaki:Lie}), and the degree~$0$ curves are
positive roots only. This partition needs a geometric insight; we
address this issue in \autoref{lem.extensible} below.
For now, merely recall (see \latin{loc.\ cit}.) that, for a root system
$\Sigma:=\root_0(N,h)$, the following objects determine one another:
\roster*
\item
a partition $\Sigma=P\cup(-P)$ into positive and negative roots,
\item
a Weyl chamber $\Weyl$ for the Weyl group of~$\Sigma$, and
\item
a linear functional $\ell\:\Z\Sigma\to\R$ such that $\ell(r)\ne0$ for each
$r\in\Sigma$ (or, more precisely, a connected component of the space of such
functionals).
\endroster

\algorithm\label{alg.Vinberg}
By \cite[\S\,8]{Huybrechts},
the \emph{smooth}
rational curves are the
(outward vectors orthogonal to the) walls of the fundamental polyhedron (of
the group generated by reflections) $\Nef X\ni h$,
and the latter is computed via Vinberg's algorithm
\cite{Vinberg:polyhedron}
\emph{provided that the partition~\eqref{eq.partition} is known}.
Denote by $\Delta_0$ the set of walls of the corresponding Weyl chamber~$\Weyl$.
Then, the smooth rational curves of degree~$d$ are the elements of the sets
$\Delta_d$ defined recursively:
\[*
\Delta_d:=\bigl\{u\in\root_d(N,h)\bigm|
 \text{$u\cdot v\ge0$ for each $v\in\Delta_r$, $r<d$}\bigr\}.
\]
We define the \emph{Fano graph} of~$X$ as
\[*
\Fn^dX:=
 \bigcup_{r=0}^d\Delta_r,
\]
where two vertices $u,v$ are connected by an edge of multiplicity
$u\cdot v\ge0$ (no edge if $u\cdot v=0$).
In other words, $\Fn^dX$ is the dual adjacency graph of the smooth rational
curves on~$X$ of degree up to~$d$.
When working with lines, we use the ``traditional'' shortcuts
\[*
\alignedat2
\Fnex X&:=\Fn^1X\:&\quad&
 \text{the \emph{extended Fano graph}},\\
\Fn X&:=\Fn^1X\sminus\Fn^0X\:&&
 \text{the \emph{\rom(plain\rom) Fano graph}, or \emph{graph of lines}}.
\endalignedat
\]
In view of \autoref{lem.extensible} below, typically we start from a
``prescribed'' graph of lines $\Fn X$, whereas $\Fnex X$ is the ultimate
result of our computation.
\endalgorithm

We extend this terminology and notation to an arbitrary polarized lattice
$N\ni h$. However, this time we need to keep track of the chosen Weyl
chamber~$\Weyl$; therefore, unless $\root_0(N,h)=\varnothing$, we speak
about
\[*
\Fn_\Weyl^d(N,h),\quad
\Fnex_\Weyl(N,h),\quad
\Fn_\Weyl(N,h),
\]
\etc.
In the most interesting case where $N=\Fano_{2n}(\graph)$ for a
graph~$\graph$,
\autoref{lem.extensible} below provides a canonical choice of~$\Weyl$; in
this case, $\Weyl$ is typically omitted. Likewise, we always omit
$\Delta=\varnothing$ when speaking about smooth polarized $K3$-surfaces.

\remark \label{rem.vocabulary}
By abuse of language we speak of a \emph{smooth graph/configuration} $\graph$ when there exists a \emph{smooth} \KK-sextic $X$, such that $\graph$ is its Fano graph. In the same manner we speak of \emph{special/non-special graph/configuration} (for the definition of special sextic see \autoref{s.special} below).
\endremark

\subsection{Projective models}\label{s.models}
A non-degenerate hyperbolic lattice $N$ equals $\NS(X)$ for some abstract
$K3$-surface~$X$ if and only if $N$ admits a primitive isometry to the
\emph{$K3$-lattice}
\[*
\bL:=H_2(X)\same2\bE_8\oplus3\bU.
\]
We call such lattices \emph{geometric}.
The existence of a primitive isometry is established \via\
\cite[Theorem 1.12.2]{Nikulin:forms} in terms of the
discriminant forms.

The following definition is context dependent: we decided not to alter the
established terminology.

\definition\label{def.admissible}
A polarized lattice $N\ni h$ is said to be \emph{$m$-admissible}, $m=1,2,3$, if
there is no vector $e\in N$ such that
\roster
\item\label{i.isotropic}
$e^2=0$ and $e\cdot h=r$ (\emph{$r$-isotropic vector}) for   $\ls|r|\le m$.
\endroster
In the context of \emph{smooth} polarized $K3$-surfaces, we require, in
addition, that there is no vector $e\in N$ such that
\roster[\lastitem]
\item\label{i.exceptional}
$e^2=-2$ and $e\cdot h=0$ (\emph{exceptional divisor}).
\endroster
(For hyperelliptic models the smoothness is understood as that of the
ramification locus.) The presence of ``bad'' vectors is easily detected by
\autoref{alg.LLL}.

A lattice $N\ni h$ is called \emph{$m$-geometric} if it is $m$-admissible and
geometric, \ie, admits a primitive isometry to~$\bL$.
\enddefinition

We refer to \cite[\S\,2.4 and Lemma~2.27]{degt.Rams:octics} (see also references
therein) for a detailed discussion of the significance of the
$m$-admissibility depending on the degree of the quasi-polarization. In the
present paper, we are concerned with $2$-admissible $6$-polarized
lattices. (One of our results is \autoref{th.special} below providing a geometric
interpretation of the $3$-admissibility in the context of sextics.)
We omit the prefix $m$- whenever it is understood from the context.

\lemma[see {\cite[Lemma~2.27]{degt.Rams:octics}}]\label{lem.geometric.sextic}
A $6$-polarized lattice $N\ni h$ is isomorphic to $\NS(X)\ni h$ for a
birational $K3$-sextic $X\to\Cp4$ if and only if it is $2$-geometric.
\done
\endlemma

As mentioned above, we are mainly interested in the lattices of the form
$\Fano_{2n}(\graph)$ or finite index extensions $N\supset\Fano_{2n}(\graph)$
thereof; in the latter case, we speak about the pair $(\graph,\CK)$
and denote $N=\Fano_{2n}(\graph,\CK)$, where the isotropic subgroup
\[*
\CK:=N/\Fano_{2n}(\graph)\subset\discr\Fano_{2n}(\graph)
 :=\Fano_{2n}(\graph)\dual\!/\Fano_{2n}(\graph)
\]
is the \emph{kernel} of the extension (see
\cite[Proposition~1.4.1]{Nikulin:forms}).
In other words, we deal with a lattice $N$ generated over~$\Q$ by~$h$ and the
vertices of~$\graph$.
We assume that $\graph$ is $\{1,\ldots,d\}$-colored and try to find a
$2n$-quasi-polarized $K3$-surface~$X$ such that
$\graph\subset\Fn^dX\sminus\Fn^0X$.
Below, we site a few definitions and statements from~\cite{degt.Rams:octics};
as in \latin{loc.\ cit}., we confine ourselves to the case $d=1$, \ie,
lines.
Thus, $\graph$ is a \emph{plain} graph, with all vertices of degree~$1$.

\definition\label{def.separating}
A Weyl chamber~$\Weyl$ for $\root_0(N,h)$ is called \emph{compatible} with
a subset $\graph\subset\root_1(N,h)$
if $\graph\subset\Fn_\Weyl(N,h)$. A root $e\in\root_0(N,h)$ is called
\emph{separating} (for~$\graph$) if there is a pair of vertices
$u,v\in\graph$ such that $u\cdot e<0<v\cdot e$.
\enddefinition

\lemma[see {\cite[Lemma~2.10]{degt.Rams:octics}}]\label{lem.extensible}
A pair $(\graph,\CK)$ admits a compatible Weyl chamber if and only if it is
\emph{extensible}, \ie, has no separating roots. In this case, a compatible
Weyl chamber is unique\rom: positive roots are those satisfying
$e\cdot\sum\graph>0$.
\done
\endlemma

\remark
The definitions above extend literally to any set
$\graph:=\graph_1\cup\ldots\cup\graph_d$, $\graph_k\subset\root_k(N,h)$.
The extensibility is necessary for the existence of a compatible Weyl
chamber, and an analogue of \autoref{lem.extensible} gives us a unique
Weyl chamber~$\Weyl$ in $N:=\Fano_{2n}(\graph,\CK)$ such that
$\graph_1\subset\Fn_\Weyl(N,h)$ and all vertices $v\in\graph_k$, $k\ge2$,
represent rational curves with all components of positive degree.
However, we can no longer guarantee that these curves are irreducible, \cf.
\cite[Example~8.4]{degt:conics}.
 For this reason, in the few cases where
conics are considered, their irreducibility is to be reconfirmed explicitly
by \autoref{alg.Vinberg}. For example, see \autoref{s.equiconic} below.
\endremark

\definition\label{def.graph}
With the degree~$2n$ fixed, a pair $(\graph,\CK)$ is called
\emph{$m$-admissible} if it is extensible and the lattice
$\Fano_{2n}(\graph,\CK)$ is $m$-admissible. If, in addition,
$\Fano_{2n}(\graph,\CK)$ is geometric, then $\graph$ and $(\graph,\CK)$ are
called \emph{$m$-subgeometric} and $\CK$ is called an \emph{$m$-geometric
kernel} for~$\graph$. Finally, $\graph$ is $m$-geometric if it is
$m$-subgeometric and
\emph{saturated}, \ie,
\[*
\graph=\sat_{2n}(\graph,\CK):=\Fn\Fano_{2n}(\graph,\CK)
\]
for some $m$-geometric kernel~$\CK$. The following statement is immediate.
\enddefinition


\theorem[see {\cite[Theorem~3.9]{degt.Rams:octics}}]\label{th.geometric}
A plain graph $\graph$ is $2$-geometric in a certain  degree $2n\ge4$
if and only if
$\graph\same\Fn X$ for some degree~$2n$ birational quasi-polarized
$K3$-surface $X\to\Cp{n+1}$
such that $\NS(X)\otimes\Q$ is generated by $h$ and the
lines.
\done
\endtheorem

A similar statement holds for $\{0,1\}$-colored graphs \vs. extended Fano
graphs $\Fnex X$ of $K3$-surfaces, \cf.
\cite[Theorem~3.10]{degt.Rams:octics}. If the $K3$-surface in
\autoref{th.geometric} is to be smooth, the extensibility condition should be
replaced with the stronger condition~\iref{i.exceptional} in the
\autoref{def.admissible} of admissibility.

\algorithm\label{alg.Nikulin}
The geometric kernels for a given graph~$\graph$ are found among the
isotropic subgroups of the finite abelian group
$\discr\Fano_{2n}(\graph)$, see
Proposition~1.4.1 in \cite{Nikulin:forms} for the relevance and
Theorem~1.12.2 in \latin{loc.\ cit}.\ for primitive isometries.
We make use of the symmetry group $\Aut\graph$ and exclude the (orbits of)
isotropic classes containing $r$-isotropic vectors, $\ls|r|\le m$, as in
\autoref{def.admissible}\iref{i.isotropic} or
\roster*
\item
separating roots as in \autoref{def.separating} (if singularities are
allowed) or
\item
exceptional divisors as in \autoref{def.admissible}\iref{i.exceptional} (in
the smooth setting),
\endroster
as well as all subgroups containing such classes.
``Bad'' vectors are detected using \autoref{alg.LLL};
many examples are found in
\autoref{S.bounds} or in
the proof of \autoref{th.ramification}.
As mentioned above, the \emph{existence} of a ``bad'' vector is proved easily:
we present one, as in the examples;
for the non-existence, we have to rely upon the software.
\endalgorithm

Thus, all conditions leading to the assertion that a given graph is
geometric are algorithmically verifiable; we refer to
\cite[\S\,B.1 and \S\,B.2]{degt.Rams:octics} for details and tweaks.

\remark\label{rem.ext}
Given a plain graph~$\graph$ and a geometric kernel~$\CK$, one can use
\autoref{lem.extensible} and \autoref{alg.Vinberg} to compute the extended
graph
\[*
\graphex:=\Fnex\Fano_{2n}(\graph,\CK).
\]
The assignment $\graph\mapsto\graphex$ is not a functor, see, for example,
$\Phi_{30}''$ \vs. $\Phi_{30}^{5\prime\prime}$ or
$\Lambda^\mathrm{A}_{24}$ \vs. $\Lambda^4_{24}$ in \autoref{lemma:easygraphs}
and \autoref{rem.Lambda}. However, using \autoref{alg.Nikulin}, we assert
that each large ($\ls|\graph|\ge36$) geometric graph~$\graph$ in degree~$6$
has a unique geometric kernel $\CK=0$ and, hence, a unique extension $\graphex$.
\endremark

\subsection{The moduli space}\label{s.moduli}
Let $X$ be a $2n$-quasi polarized $K3$-surface. Fix its
Fano graph $\graph:=\Fn^dX$ and let
\[*
N := (\Q\graph+\Q h)\cap\NS(X)\same\Fano_{2n}(\graph,\CK)
\]
for an appropriate kernel~$\CK$.
It is immediate from Vinberg's \autoref{alg.Vinberg} that the
equilinear (equiconical, \etc.) stratum
\[
\bigl\{\text{$Y\to\Cp{n+1}$ a $2n$-quasi-polarized $K3$-surface}\bigm|
 \Fn^dY\same\graph\bigr\}
\label{eq.stratum}
\]
containing~$X$
is an open subset of the respective stratum of the space of lattice
$N$-polarized $K3$-surfaces in the sense of
Dolgachev~\cite{Dolgachev:polarized}. Thus, all strata~\eqref{eq.stratum} are
in a bijection with the following sets of data:
\roster*
\item
an $(\Aut\graph)$-orbit of geometric lattices $N:=\Fano_{2n}(\graph,\CK)$
(or kernels $\CK$) with the property that $\Fn^d(N,h)=\graph$, \cf.
Algorithms~\ref{alg.Nikulin} and~\ref{alg.Vinberg}, and
\item
for each~$N$, a primitive isometry $N\to\bL$ up to the left-right action of
the group $\OG_{h}(\NS)\times\OG^+(\bL)$, where
$\OG^+$ is the subgroup preserving the
\emph{positive sign structure}, \ie, coherent orientation of maximal
positive subspaces.
\endroster

For the latter, we use Nikulin's theory~\cite{Nikulin:forms}: the genus of
the \emph{transcendental lattice} $T:=N^\perp\subset\bL$ is determined by~$N$
and, with~$T$ fixed, primitive isometries $N\to\bL$ are in a bijection with
the $\OG_h(N)\times\OG^+(T)$-orbits of anti-isomorphisms
\[
\phi\:\discr N\longto\discr T,
\label{eq.phi}
\]
see \cite[Proposition~1.6.1]{Nikulin:forms}.
One way to define~$\phi$ is that, for $u\in N\dual$ and $v\in T\dual$,
\[
\phi(u\bmod N)=v\bmod T\quad\text{if and only if}\quad u+v\in\bL.
\label{eq.phi.def}
\]

\algorithm\label{alg.aut}
The finite group $\OG_h(\Fano_{2n}(\graph))=\Aut\graph$ is easily computed using the
\texttt{digraph} package in \GAP~\cite{GAP4.13}. Then,
$\OG_h(\Fano_{2n}(\graph,\CK))$ is the stabilizer of~$\CK$ under the
canonical action of $\Aut\graph$ on $\discr\Fano_{2n}(\graph)$.

Both $\operatorname{genus}(T)$ and
$\OG(T)$ are a classical subject (Gauss~\cite{Gauss:Disquisitiones})
if $T$ is positive definite of rank~$2$
(\ie, $X$ is a singular $K3$-surface). If $\rank T\ge3$, we use
Miranda--Morrison
theory~\cite{Miranda.Morrison:1,Miranda.Morrison:2,Miranda.Morrison:book},
which combines
$\Coker\bigl[\OG^+(T)\to\Aut(\discr T)\bigr]$
and $\operatorname{genus}(T)$
into a single finite
$2$-elementary abelian group determined by~$N$.
\endalgorithm


\subsection{Projective automorphisms}\label{s.aut}
Let $X$ be a $2n$-quasi-polarized $K3$-surface. We denote $N:=\NS(X)$ and let
$T:=N^\perp\subset\bL$ be the transcendental lattice.

If $X$ is singular ($\rank N=20$), its group
$\Aut_hX$ of projective automorphisms is
\[*
\Aut_hX=\OG_h(N)\times_\phi\OG^+(T),
\]
see~\eqref{eq.phi},~\eqref{eq.phi.def} for the definition of~$\phi$ and
\autoref{alg.aut} for the
computation.
Otherwise, $\rank N\le19$, one has
\[*
\Aut_hY=\OG_h(N)\times_\phi\{\pm\id_T\}
\]
for a \emph{very general} representative~$Y$ of the respective
lattice $N$-polarized stratum, and it is these groups that are stated in
\autoref{tab.list}.
If $\rank N$ is odd, instead of the vague ``very general''
it suffices to assume that $\NS(Y)=N$.

\subsection{The line-by-line algorithm}\label{s.line.by.line}
In many proofs, we start with a sufficiently large graph~$\graph$ and
extend it by adding several new vertices, one extra vertex~$u$ at a time.
We assume known (from the geometry of the problem) that $u\cdot v\in\{0,1\}$
for each $v\in\graph$. Then, the new lattice $\Fano(\graph\sqcup v)$
is determined by the \emph{support}
\[
\|u\|:=\bigl\{v\in\graph\bigm|u\cdot v=1\bigr\}\subset\graph
\label{eq.supp}
\]
of~$u$, which is merely a subset of~$\graph$. If several vertices
$u_1,u_2,\ldots$ are to be added at once, in addition to the multiset
$\bold{u}:=\{\|u_1\|,\|u_2\|,\ldots\}$ of subsets of~$\graph$ we need to specify the
Gram matrix
\[
\fm:=[u_i\cdot u_j].
\label{eq.fm}
\]
For the new graph
$\graph':=\graph\sqcup\bold{u}(\fm)$,
we compute the new
lattice $N':=\Fano_{2n}(\graph')$ and analyse its properties, most notably
$\Gs_+N=1$, admissibility, extensibility, \etc.

We refer to \cite[\S\,B3 and \S\,B4]{degt.Rams:octics} for the details
concerning the implementation; a plethora of examples is found in
the other appendices in
\latin{loc.\ cit}. Certainly, we make full use of the group $\Aut\graph$ of
symmetries and of geometric insight limiting the values taken by the
supports~\eqref{eq.supp} and matrices~\eqref{eq.fm}; in the present paper, a
good example is the proof of \autoref{th.Humbert}.


\subsection{The Kummer index of a quasi-polarized $K3$-surface}\label{s.Kummer}
We conclude this section by proving a general statement announced in
\cite[\S\,6.4]{degt.Rams:octics}. It can be see as an analogue of
\cite[Corollary 2]{Nikulin:Kummer} and \cite{Barth:classninecusps}.

\theorem\label{th.Kummer}
A quasi-polarized $K3$-surface $(X,h)$ with $h^2=2\bmod4$ contains at most
$12$ pairwise disjoint
smooth rational curves of odd degree.
\endtheorem

\proof
Assume that there are $13$ curves $l_i\in\NS(X)$ of odd degrees $d_i$,
$i=1,\ldots,13$, and consider
the sublattice $S\subset\NS(X)$ spanned by these curves. The
discriminant group
$\discr S\same(\Z/2)^{13}$ is generated by the pairwise orthogonal classes
\[*
\Ga_i:=\frac12l_i\bmod S,\quad
\Ga_i^2=-\frac12\bmod2\Z,\quad
i=1,\ldots,13.
\]
Up to automorphisms (essentially, ignoring the reflections, up to $\SG{13}$),
there are three classes of isotropic vectors, \viz. $\Ga_1+\ldots+\Ga_n$,
$n=4,8,12$, and the former contains
a root, \eg, $e:=\frac12(l_1+\ldots+l_4)$, which is positive whenever
all $l_i$ are. The vectors $e,l_1,\ldots,l_4$ span~$\bD_4$,
implying that, in the
finite index extension by~$e$, all four curves $l_1,l_2,l_3,l_4$ cannot be walls of a
common Weyl chamber (\cf. Nikulin~\cite{Nikulin:Kummer} and \autoref{alg.Vinberg}).

We conclude that the only finite
index extension $\tilde{S}\supset S$ admitting a primitive isometry to~$\bL$
and having all $13$ curves $l_1,\ldots,l_{13}$
among the walls of a common Weyl chamber is that by,
say,
\[*
\kappa_1:=\Ga_1+\ldots+\Ga_8,\qquad\kappa_2:=\Ga_5+\ldots+\Ga_{12}.
\]
Then
\[
T_S:=\tilde{S}^\perp\same6\bA_1\oplus3\bA_1(-1)\subset\bL
\label{eq.TS}
\]
is the double of an odd unimodular lattice.
Now, $h$ projects to the vector
\[*
\tilde{h}:=h+\frac12\sum_{i=1}^{13}d_il_i\in T_S\dual.
\]
The difference $\tilde{h}-h$ is in $S\dual$ and its
class $\Gd=\sum_{i=1}^{13}\Ga_i\bmod S\in\discr S$ survives to
$\discr\tilde{S}$, as $\Gd\cdot\kappa_1=\Gd\cdot\kappa_2=0\bmod\Z$.
Furthermore, $\Gd$ is characteristic in $\discr S$, and so
it is in $\discr\tilde{S}$. Hence, the class $\phi(\Gd)$ of~$\tilde{h}$, \cf.
\eqref{eq.phi} and~\eqref{eq.phi.def}, is characteristic in $\discr T_S$ and,
considering~\eqref{eq.TS},
$2\tilde{h}$ is a characteristic vector in the unimodular lattice
$T_S(\frac12)$. We conclude that
\[*
\frac12(2\tilde{h})^2=\sum_{i=1}^{13}d_i^2+2h^2\equiv13+2h^2\equiv\Gs(T_S)=-3\bmod8
\]
(recall that $d_i^2=1\bmod8$)
and, hence, $h^2=0\bmod4$.
\endproof

\remark
(1) Apart from being of interest on its own, \autoref{th.Kummer}
plays a certain r\^{o}le in our
arguments:
a configuration
with at least $13$ disjoint lines
is discarded
immediately, without waiting for the late stages of the algorithm.

\smallskip
(2) In general, \autoref{th.Kummer} provides a useful insight that also explains certain aspects of classification of large line configurations for other polarizations $h$ such that  $h^2=0\bmod4$ (\cf.~the general approach presented in \cite{degt.Rams:octics}).
\endremark


\section{Simple bounds for sextics}\label{S.bounds}

In the remainder of this paper
(except \autoref{s.sextics})
	 we restrict our attention to \KK-sextics, so
that all polarized lattices
considered are  $6$-polarized.
We omit the degree of the polarization when it leads to no ambiguity;
\eg, we put~$\Fano(\graph):= \Fano_{6}(\graph)$.

Recall that
a
\emph{positive root}
in the Neron-Severi lattice $\NS(X)$
is an exceptional
divisor, \ie, a union of $(-2)$-curves that are contracted by
the map $X\to\Cp4$,
see,
\eg, \autoref{s.LLL} or \cite[\S\,2.3]{degt.Rams:octics}.
Positive roots are detected by \autoref{lem.extensible}.
Similarly, a \emph{degree~$1$ curve} is an element of
$\root_1\NS(X)$:
it has the form
$l+e_1+e_2+\ldots$, where $l$ is a line and $e_i$ are  pairwise
disjoint positive roots.
By \cite[Lemma~2.25]{degt.Rams:octics}, we have $l\cdot e_i=1$.
Below, we use the terms ``positive root'' and ``degree~$1$ curve'' for,
respectively,
``exceptional divisors'' and ``lines'' about which we cannot immediately
assert the irreducibility.

\newcommand{\Xh}{X}

\subsection{Stars of vertices}\label{s.stars}
Here we collect various constraints on the number of edges of
the (extended) Fano
graph of a \KK-sextic that share a vertex.
We define
the \emph{star} of a vertex  $v \in \graph$
and its \emph{valency} \via
$$
\Star v := \bigl\{l \in \graph  \bigm| l \cdot v = 1\bigr\} \quad \mbox{and}
\quad \val v := \ls|\Star v| \, .
$$

\proposition\label{prop.exceptional}
In a $K3$-sextic $\Xh$,
an exceptional divisor may intersect at most six
\rom(necessarily pairwise disjoint\rom) lines.
\endproposition

\proof
The statement follows from the Hodge index theorem: should there be seven
lines $l_i$, the lattice spanned by $h$, $e$, and $l_1,\ldots,l_7$ would have
$\Gs_+=2$.
\endproof

\remark\label{rem.exceptional}
If an exceptional divisor~$e$ intersects five lines $l_1,\ldots,l_5$,
then there is a
sixth degree~$1$ curve
\[*
l_6=h-3e-(l_1+\ldots+l_5)
\]
that also intersects~$e$.
\endremark

\theorem\label{th.star}
Let $\graph:=\Fn(\Xh)$ be the plain Fano graph of a $K3$-sextic~$\Xh$.
Then the star of any vertex $v\in\graph$ is of one of the following three
types\rom:
\roster
\item\label{star.0}
$l_1,\ldots,l_k$ \rom(\ie, $k\bA_1$\rom) with $0\le k\le9$, or
\item\label{star.1}
$m_1,m_2,l_1,\ldots,l_k$ \rom(\ie, $\bA_2\oplus k\bA_1$\rom)
with $0\le k\le9$, $k\ne8$, or
\item\label{star.2}
$m_1,m_2,n_1,n_2$ \rom(\ie, $2\bA_2$\rom),
\endroster
where
\[*
m_1\cdot m_2=n_1\cdot n_2=1,\quad
m_r\cdot n_s=m_r\cdot l_i=l_i\cdot l_j=0
\]
for $r,s=1,2$ and $i,j=1\ldots,k$, $i\ne j$.
Furthermore,
\roster[\lastitem]
\item\label{star.1.more}
in case~\iref{star.1} with $k=9$,
one has $3h=5v+2(m_1+m_2)+l_1+\ldots+l_9$\rom;
\item\label{star.2.more}
in case~\iref{star.2}, $e:=h-2v-m_1-m_2-n_1-n_2$ is
a positive root.
\endroster
\endtheorem

\remark\label{rem.star}
In case~\iref{star.2} of \autoref{th.star},
the positive root in item~\iref{star.2.more} implies that the sextic has
singular points
and that no other line
intersects either~$v$ or more than one of $m_1,m_2,n_1,n_2$.
\endremark

\proof[Proof of \autoref{th.star}]
Statements~\iref{star.1.more},~\iref{star.2.more} are immediate, see
\autoref{obs.ker}.

Consider the graph~$\graph$ consisting of a vertex~$v$ and ten pairwise
disjoint vertices $l_1,\ldots,l_{10}$ such that $l_i\cdot v=1$ and let
$N:=\Fano(\graph)$.
The group $\discr_2N\same(\Z/2)^{10}$ is generated by the ten
isotropic vectors
\[*
\Gl_i:=\frac12h+\frac12v+\frac12\sum_{s\ne i}l_s,\quad i=1,\ldots,10,
\]
so that $\Gl_1\cdot\Gl_j=\frac12\bmod\Z$ for $i\ne j$. Up to
isomorphism, the isotropic classes are $\Gl_1+\ldots+\Gl_r$,
$r=1,4,5,8,9$, most of which contain forbidden vectors:
\begin{alignat*}2
r=4\:&\quad e=\tfrac12(l_1-l_2+l_3-l_4)&\quad&\text{(separates $l_1$ and $l_2$)},\\
r=5\:&\quad e=\tfrac12(h-v-l_6-l_7-\ldots-l_{10})&&\text{(separates $v$ and $l_{10}$)},\\
r=8\:&\quad e=h-2v-\tfrac12(l_1+\ldots+l_{8})&&\text{(separates $v$ and $l_{10}$)},\\
r=9\:&\quad p=\tfrac12(h-v-l_{10})&&\text{($2$-isotropic vector)}.
\end{alignat*}
Since $N$ itself does not admit a primitive isometry into~$\bL$
(as $\det N=3\cdot2^{10}$), the only geometric extension
\smash{$\tilde{N}\supset N$}
is that
by, say, $\Gl_{10}$; upon renaming $l_{10}$ to $m_1$, it has an extra
degree~$1$ curve~$m_2$ given by the relation in item~\iref{star.1.more}.

We still need to show that the linear component $m_2'$ of~$m_2$ intersects
both~$v$ and $l_{10}=m_1$, which is distinguished by our choice of the
extension $\tilde{N}\supset N$. To this end, we extend $\graph$ and
$N,\tilde{N}$ by an extra vertex/vector~$m_2'$.
We need that $m_2'\cdot\Gl_{10}\in\Z$; hence, $m_2'\cdot v=1$ and, assuming
that $m_2'\cdot l_{10}=0$, the new graph $\graph':=\graph\cup m_2'$ and
lattice $N':=\Fano(\graph')=N+\Z m_2'$ are as in item~\iref{star.0} with
$k=11$. This time, the group $\discr_2N'=(\Z/2)^{10}\oplus\Z/4$ is generated
by the same ten vectors $\Gl_1,\ldots,\Gl_{10}$ and
$\eta:=\frac34(h+l_1+\ldots+l_{10})+\frac14m_2'$. It has the same five orbits
of isotropic classes, containing the same forbidden vectors. Hence, the
lattice \smash{$\tilde{N}':=\tilde{N}+\Z m_2'$} has no proper admissible finite index
extensions. This lattice is not geometric.

This argument shows that $k\le9$ in cases~\iref{star.0} and~\iref{star.1}.
Assume that $k=8$ in case~\iref{star.1} and let $l_9$ be the degree~$1$ curve
given by the relation in item~\iref{star.1.more}. If we assume that the
linear component $l_9'$ of~$l_9$ is disjoint from~$v$, adjoining the
vertex~$l'_9$ to the graph makes it non-hyperbolic.

Formally, in case~\iref{star.1} we also need to rule out an extra line~$n$
such that either
\roster*
\item
$n\cdot\ell=n\cdot m_1=n\cdot m_2=1$ (the graph is not hyperbolic) or
\item
$n\cdot\ell=n\cdot m_1=1$, $n\cdot m_2=0$ ($e:=h-\ell-m_1-m_2-n$
is $2$-isotropic).
\endroster

Finally, in case~\iref{star.2}, there are no other lines in $\Star v$ due to
\autoref{rem.star}.
\endproof

\remark\label{rem.conics}
Although we do not systematically study conics in this paper, in
\autoref{S.aut} below we need to consider conics that intersect a certain
line at two points. We observe that all conclusions of \autoref{th.star} and
\autoref{rem.star} hold literally if we replace the ``triangles'' $m_1,m_2$
and $n_1,n_2$ with conics $c:=m_1+m_2$, $d:=n_1+n_2$, so that
$c\cdot h=c\cdot v=d\cdot h=c\cdot v=2$.
\endremark


\subsection{Biquadrangles}\label{s.biquads}

Recall that a biquadrangle
is a complete bipartite subgraph $K_{2,3}$ (see also \autoref{rem:defbiquadrangles} below). 
By \cite[Lemma~5.2]{degt.Rams:octics},
the Fano graph of a \KK-octic can contain such a subgraph only if
the surface in question is special. 
 In contrast, biquadrangles play important role in our study
of Fano graphs of \KK-sextics.
In particular, they provide a useful constraint on the graphs,
see \autoref{rem.prevalidate} below.

\theorem\label{th.biquad}
Let $\graph:=\Fn(\Xh)$ be the plain Fano graph of a $K3$-sextic~$\Xh$.
Then, for two disjoint vertices $l_1,l_2\in\graph$, the intersection
$I:={\Star l_1}\cap{\Star l_2}$ of their stars \rom(which is
discrete due to \autoref{th.star}\rom) consists of at most three vertices.
If $I=\{m_1,m_2,m_3\}$ does have three vertices, then
\[*
\bl_3:=h-(l_1+l_2+m_1+m_2+m_3)
\]
is a degree~$1$ curve
such that
\[*
\bl_3\cdot l_1=\bl_3\cdot l_2=0\quad\text{and}\quad
\bl_3\cdot m_i=1\quad\text{for $i=1,2,3$}.
\]
\endtheorem

\remark \label{rem:defbiquadrangles}
A configuration $l_1,l_2,m_1,m_2,m_3$ as in \autoref{th.biquad}, \ie, such
that
\[\label{eq.biquad}
l_i\cdot l_j=m_r\cdot m_s=0\quad\text{and}\quad
l_i\cdot m_r=1\quad\text{for $i\ne j$ and $r\ne s$},
\]
where $1\le i,j\le2$, $1\le r,s\le3$,
is called a \emph{biquadrangle}.
In a smooth configuration, a biquadrangle is ``completed'' by the sixth
line $l_3:=\bl_3$ as in \autoref{th.biquad}, disjoint from~$l_1$, $l_2$ and
intersecting all~$m_r$, $r=1,2,3$. In general, the degree~$1$ component~$l_3$
of~$\bl_3$ may be also disjoint from some of~$m_i$.
\endremark

\proof[Proof of \autoref{th.biquad}]
If $I\supset\{m_1,\ldots,m_4\}$, then the root
\[*
h-(l_1+l_2+m_1+\ldots+m_4)
\]
separates $l_1,l_2$ from $m_1,\ldots,m_4$.
The other statements are immediate.
\endproof

A geometric interpretation of \autoref{cor.biquad} below can be found
in \cite[\S\,2.3]{DDK}.

\corollary\label{cor.biquad}
Assuming that
$l_1,l_2,m_1,m_2,m_3$ is a biquadrangle
and denoting by~$l_3$
the linear component
of the degree~$1$ curve $\bl_3$ as in \autoref{th.biquad}, we
conclude that any other line intersects at most one of
the six lines $l_i,m_i$, $1\le i\le3$. If $\bl_3=l_3$ is irreducible,
any other line intersects \emph{exactly} one of $l_i,m_i$.
\done
\endcorollary

The case $\bl_3=l_3$ of \autoref{cor.biquad}, \ie, a collection
\[
\biquad:=\{l_1,l_2,l_3,m_1,m_2,m_3\}\subset\graph
\label{eq.complete}
\]
of six lines intersecting as in~\eqref{eq.biquad} is called a
\emph{complete biquadrangle}. \autoref{cor.biquad} states that each line
$n\in\graph\sminus\biquad$ intersects \emph{exactly one} element
of~$\biquad$. If $X$ is smooth, any biquadrangle~\eqref{eq.biquad} is part of
a unique complete biquadrangle~\eqref{eq.complete}.

\remark\label{rem.prevalidate}
Since separating roots would appear at a relatively late stage of our computation,
to save time we use \autoref{th.biquad} and other similar statements (\cf.
\autoref{lem.trig.trig} or \autoref{lem.trig.sec} below) to pre-validate
(the Gram matrix of) a (partial) graph: those with a forbidden subgraph are
rejected immediately.
\endremark


\subsection{Taxonomy of graphs and pencils}\label{s,pencils}

As in \cite{degt:lines,degt.Rams:octics},
given a graph $\graph$, we
consider
the lattice $\ZZ\graph$,
\cf.~\eqref{eq.Fano}, and call
$\graph$ {\emph hyperbolic} (resp.\ \emph{parabolic}, resp.\ \emph{elliptic})
if $\sigma_{+}(\ZZ\graph) = 1$
(resp.\  $\sigma_{+}(\ZZ\graph) = 0$ and  $\sigma_{0}(\ZZ\graph) > 0$,
 resp.\  $\sigma_{+}(\ZZ\graph) = \sigma_{0}(\ZZ\graph) = 0$). For
 a parabolic or elliptic graph $\graph$ we put
 $\mu(\graph):= \rank(\ZZ\graph/\ker)$ to denote its Milnor number and introduce  an order on  the set of
 (isomorphism classes of)
 connected parabolic graphs (\latin{a.k.a}.\ affine Dynkin diagrams):
if $\mu:=\mu(\Sigma') =\mu(\Sigma'')$, we define
$\tA_{\mu} < \tD_{\mu} < \tE_{\mu}$; otherwise, $\Sigma'<\Sigma''$ whenever
$\mu(\Sigma') <  \mu(\Sigma'')$.


A \emph{minimal fiber} $\fiber$ of a hyperbolic
or parabolic graph~$\graph$ (see \cite[Definition~4.1]{degt.Rams:octics}) is a  connected parabolic
subgraph of $\graph$ that is minimal with respect to the order ``$<$'' introduced above.
We put
$\kappal_\fiber := \sum_{w \in \fiber} n_w w$, $n_v>0$,
to denote the \emph{fundamental cycle}
of~$\Sigma$, see \cite[(4.1)]{degt.Rams:octics}.
Besides,
as in \cite{degt:lines,degt.Rams:octics}, for
a
 hyperbolic graph  $\graph$ with a connected parabolic
subgraph $\fiber$, we consider
the
\emph{pencil containing $\fiber$}, \ie,
the maximal parabolic subgraph
of  $\graph$ that contains $\fiber$:
\[*
\pencil:=\fiber\cup\bigl\{v\in\graph\bigm|
\text{$v\cdot l=0$ for all $l\in\fiber$}\bigr\}
\]
and the sets
\begin{alignat*}2
\operatorname{Sec}\fiber&
 :=\bigl\{v\in\graph\sminus\fiber\bigm|v\cdot {\kappal_\fiber} > 0\bigr\}
 &\quad&\text{(\emph{sections} of~$\pencil$)},\\
\sec\fiber&
 :=  \bigl\{v\in\graph\sminus\fiber\bigm|v\cdot {\kappal_\fiber} =1 \bigr\}
 &\quad&\text{(\emph{simple sections})},\rlap{ and}\\
\sec e&
 :=\bigl\{v\in\graph\sminus\fiber\bigm|v\cdot e={v\cdot\kappal_\fiber}=1\bigr\} 
 &\quad&\text{(for $e\in\fiber$)}.
\end{alignat*}
The \emph{multiplicity} of a section
$v\in\operatorname{Sec}\fiber$ is $v\cdot\kappal_\fiber$, so that
simple sections are those of multiplicity~$1$.
Parabolic (connected) components of $\pencil$ are called its \emph{fibers}.
Usually, we number the vertices $\fiber=\{e_1,e_2,\ldots\}$ and abbreviate
$\sec_i:=\sec e_i$. We define
\[*
\textstyle
{\sec_i}*{\sec_j} := \max\bigl\{v\cdot\sum\sec_j\bigm|v\in\sec_i\bigr\}.
\]

\remark \label{rem:pencil.notation}
In what follows, when we speak of a graph $\graph$ of type $\fiber$
and write
$$
\graph\supset\pencil\supset\fiber \, ,
$$
we mean that $\fiber$ is a (fixed) minimal fiber and $\pencil$ is the pencil
containing~$\fiber$.
\endremark


\section{Quadrangular graphs}\label{S.quad}

In this section, we fix a quadrangular graph (\ie, a graph of the
type $\tA_3$)
\[\label{eq.quad}
\graph\supset\pencil\supset\fiber:=\{e_1,e_2,e_3,e_4\}\same\tA_3,
\]
where the vertices of~$\fiber$ are numbered cyclically. Observe that
in accordance with our taxonomy, we
assume that $\graph$ is \emph{triangular free}; in particular,
\roster*
\item
each subgraph $\sec_i$, $i=1,\ldots,4$, is discrete and
\item
there is no multi-section~$l$
\emph{with $l\cdot e_i=l\cdot e_j=1$ for $j=i\pm1\bmod4$.}
\endroster
Besides, we assume the vertices of~$\fiber$ ordered so that
\[\label{eq.quad.order}
\ls|\sec_1|\ge\ls|\sec_2|\ge\ls|\sec_4|,\quad
\ls|\sec_1|\ge\ls|\sec_3|,\quad
\ls|\sec_3|\ge\ls|\sec_4|\ \text{if $\ls|\sec_1|=\ls|\sec_2|$}.
\]

\corollary[of \autoref{th.star}]\label{cor.quad.valency}
One has $\val v\le9$ for each vertex $v\in\graph$. Hence, one also has
$\ls|\sec_i|\le7$ for each $i=1,\ldots,4$.
\done
\endcorollary

\corollary[of \autoref{th.biquad}]\label{cor.quad}
One has
\[*
{\sec_i}*{\sec_j}\le\begin{cases}
 2,&\text{if $i=j\pm1\bmod4$},\\
 3,&\text{if $\ls|i-j|=2$}.
\end{cases}
\def\qedsymbol{\donesymbol}\pushQED{\qed}\qedhere
\]
\endcorollary

\subsection{Biquadrangle free graphs}\label{s.biquad.free}
A quadrangular graph~$\graph$ is
called \emph{biquadrangle free} if it
does not contain a single biquadrangle; in other words, if the intersection
of the stars of any two disjoint vertices has cardinality at most~$2$.
In this case, the bounds of \autoref{cor.quad} automatically reduce down to
\[ \label{eq:secisecjbfree}
{\sec_i}*{\sec_j}\le\begin{cases}
 1,&\text{if $i=j\pm1\bmod2$},\\
 2,&\text{if $\ls|i-j|=2$}.
\end{cases}
\]

\lemma \label{lem:tA3biquadranglefree}
For a biquadrangle free geometric graph~$\graph$ one has either
$\graphex\same\Theta^\circ_{36}$ \rom(see \autoref{rem.ext} and
\autoref{tab.list}\rom)
or $\ls|\graph|\le35$.
\endlemma

\proof We follow  the general strategy depicted in \cite[\S\,4.4]{degt.Rams:octics}.
We start from a quadrangular biquadrangle free pencil $\pencil$
with $\ls|\pencil|\ge17$. Direct computation with the help of \cite{GAP4.13}
(essentially, we use the line-by-line algorithm and all the constraints
on the intersection of sections, \eg, \eqref{eq:secisecjbfree};
see \cite[Appendix~D]{degt.Rams:octics} for further details)
shows that  $\ls|\graph|\le35$
unless $\graph\same\Theta^\circ_{36}$.

Then, we start from a
single quadrangle $\fiber\same\tA_3$ and use the same constraints to
analyze the sets of sections
(see \cite[Appendix~C]{degt.Rams:octics}), \emph{aiming at
$\ls|\sec\fiber|\ge20$}. Once again, we arrive at a single graph
$\graph\same\Theta^\circ_{36}$. Summarizing, apart from the case
$\graph\same\Theta^\circ_{36}$, we can assume that $\ls|\pencil|\le16$ and
$\ls|\sec\fiber|\le19$; hence, $\ls|\graph|\le35$.
\endproof

\remark\label{rem.thresholds}
Note that in this or subsequent proofs we do \emph{not} assert that always
$\ls|\pencil|\le16$ or $\ls|\sec\fiber|\le19$. Rather, we find \emph{all}
graphs violating one of these bounds and, with the exception stated, we have
$\ls|\graph|\le35$ for all such graphs. \latin{Cf}.\
\cite[(4.10)]{degt.Rams:octics}.
\endremark

\subsection{Biquadrangular pencils}\label{s.biquad}
We temporarily rename the vertices of the distinguished fiber~$\fiber$,
considering the ``horizontal'' lines $l_1:=e_1$, $l_2:=e_3$ and ``vertical''
lines $m_1:=e_2$, $m_2:=e_4$, so that \eqref{eq.biquad} holds.

\definition\label{def.biquad}
A quadrangular
pencil~$\pencil$ is called \emph{biquadrangular} if it has a
\emph{bisection}, \ie, a
line~$l_0$ or~$m_0$ such that \eqref{eq.biquad} still holds.
We consider the following  four  
 types of quadrangular pencils
\roster[0]
\item\label{biquad.0}
there is neither~$l_0$ nor~$m_0$, so that the pencil is
biquadrangle free,
\item\label{biquad.1}
there is exactly one extra line $l_0$ or $m_0$,
\item\label{biquad.110}
there are both $l_0$ and $m_0$ and $l_0\cdot m_0=0$, or
\item\label{biquad.111}
there are both $l_0$ and $m_0$ and $l_0\cdot m_0=1$, so that $\fiber$ and
$l_0,m_0$ constitute a complete biquadrangle, see~\eqref{eq.complete},
\endroster
which are ordered \via\
$\iref{biquad.0}<\iref{biquad.1}<\iref{biquad.110}<\iref{biquad.111}$
(see \autoref{conv.biquad} below).
Due to \autoref{th.biquad}, the pencil cannot have more than one
horizontal bisection~$l_i$ or more than one vertical bisection~$m_r$.
In a  configuration of lines on  a smooth \KK-sextics, all biquadrangular pencils are of
type~\iref{biquad.111}, see \autoref{rem:defbiquadrangles}.
\enddefinition

\convention\label{conv.biquad}
We always assume that the quadrangular pencil~$\pencil$ in~\eqref{eq.biquad} is
\roster
\item\label{bq.max.type}
of the maximal type (among those present in~$\graph$),
\item\label{bq.max.card}
of the maximal cardinality among all pencils in~$\graph$
\emph{of the same type}, and
\item\label{bq.min}
minimal, with respect to  a certain lexicographic order,
see \eqref{eq.biquad.bound} below,
among all pencils
in~$\graph$ of the same type and cardinality.
\endroster
Convention~\iref{bq.max.type} lets us assume that all sections of~$\pencil$
are disjoint from its bisections, whereas Convention~\iref{bq.min}
is a mere technicality reducing
the overcounting.
Pencils of
type~\eqref{biquad.0} are considered in biquadrangle free graphs, see
\autoref{s.biquad.free}.
\endconvention

With~\eqref{eq.quad.order} taken into account,
\autoref{conv.biquad}\iref{bq.max.card} has the following implication for
types~\iref{biquad.1} and~\iref{biquad.111}
(but not for type~\iref{biquad.110}): whenever $l_0$ or $m_0$ is present,
one has
\[\label{eq.biquad.bound}
\val l_0\ge\val l_1\ge\val l_2,\quad
\val m_0\ge\val m_1\ge\val m_2.
\]
Indeed, if, say, $\val l_0<\val l_1$, then the fiber
$\fiber':=\{l_0,m_1,l_2,m_2\}$ would define
a pencil $\pencil'\subset\graph$ of the same type and strictly larger
cardinality (as $\pencil'$ would have fewer sections than~$\pencil$).

\remark\label{rem.biquad.bound}
Inequality~\eqref{eq.biquad.bound} is indispensable when dealing with small
pencils~$\pencil$, which require most effort. For example, if $\pencil$ is of
type~\iref{biquad.111} and $\ls|\pencil|=14$ (the worst case), it is not
uncommon that the bounds on the sequence $(\ls|\sec_i|)$, $i=1,\ldots,4$,
given by~\eqref{eq.biquad.bound} are $(4,6,4,6)$ or $(6,4,6,4)$. Since we
need at least $10$ sections and with~\eqref{eq.quad.order} taken into
account, in the former case $\pencil$ can be ignored immediately, whereas in
the latter case we are forced to confine ourselves to the relatively few
graphs with $\ls|\sec_1|=6$.
\endremark

\lemma \label{lem:tA3type1and2}
A quadrangular geometric graph~$\graph$ of type ~\iref{biquad.1}
or~\iref{biquad.110} has $\ls|\graph|\le35$.
\endlemma

\proof
Again, we follow the approach of~\cite{degt.Rams:octics} and use the obvious
identity
\[*
\ls|\graph|=\ls|\pencil|+\ls|\sec\fiber|+\Gd,
\]
where $\Gd$ is the number of bisections. We find all graphs
$\graph\supset\pencil\supset\fiber$ that are as in the statement and satisfy
either
\[*
\ls|\pencil|> M_\pencil\ \text{(\cf. \cite[Appendix~D]{degt.Rams:octics})}
\quad\text{or}\quad
\ls|\sec\fiber|> M_\fiber\ \text{(\cf. \cite[Appendix~C]{degt.Rams:octics})},
\]
where we \emph{choose}
\roster*
\item
$M_\pencil:=15$ and $M_\fiber:=19$ for $\pencil$ of type~\iref{biquad.1}
(one has $\Gd=1$) or
\item
$M_\pencil:=13$ and $M_\fiber:=20$ for $\pencil$ of type~\iref{biquad.110}
(one has $\Gd=2$).
\endroster
We find that $\ls|\graph|\le34$ for all graphs $\graph$ obtained.
Therefore, we can \emph{assume} that $\ls|\pencil|\le M_\pencil$ and
$\ls|\sec\fiber|\le M_\fiber$
(\cf. \autoref{rem.thresholds}), resulting in $\ls|\graph|\le35$.
\endproof


\lemma\label{lem.111}
If $\graph$ is a quadrangular geometric graph of type~\iref{biquad.111},
then either
\[*
\graphex\same\Theta_{36}\quad \text{\rm or}\quad \gf{N}{36},\ N=1,\ldots,15
\]
\rom(see \autoref{rem.ext} and \autoref{tab.list}\rom), or
$\ls|\graph|\le35$.
\endlemma

\proof
By \autoref{cor.biquad}, we have
\[*
\ls|\pencil|=\val l_0+\val m_0+4;\quad\text{hence},\quad
\ls|\graph|=\ls|\pencil|+\ls|\sec\fiber|+2\le3\ls|\pencil|-6
\]
by~\eqref{eq.biquad.bound} and it suffices to consider the case
$\ls|\pencil|\ge14$. This is done as explained in
\cite[Appendix~C]{degt.Rams:octics} (with \autoref{conv.biquad} and
\autoref{rem.biquad.bound} playing an important r\^{o}le), and this is the most
time/resource intensive part of our proof, resulting in the vast majority of
large graphs.
\endproof


\section{Triangular graphs and graphs of the remaining types} \label{S.other}

In this section we study
the remaining graphs, \ie, those of type other than~$\tA_3$.

\subsection{Triangular graphs}\label{s.trig}

Fix a triangular graph (\ie, one
of type $\tA_2$)
\[\label{eq.trig}
\graph\supset\pencil\supset\fiber:=\{e_1,e_2,e_3\}\same\tA_2.
\]

\corollary[of \autoref{th.star}]\label{cor.trig}
The pencil~$\pencil$ has no multisections. Furthermore,
for each $i=1,2,3$, one has either
\roster
\item\label{trig.trig}
$\sec_i=\{m_1,m_2\}$, $m_1\cdot m_2=1$, or
\item\label{trig.discrete}
$\sec_i$ is discrete and $\ls|\sec_i|\le9$.
\endroster
In the former case, $m_1,m_2$ are disjoint from all other sections
of~$\pencil$.
\done
\endcorollary

\lemma\label{lem.trig.trig}
If $\sec_i=\{m_1,m_2\}$ with $m_1\cdot m_2=1$,
\ie, case~\iref{trig.trig} of \autoref{cor.trig},
then
$\ls|\sec_i|\le6$ for $i=2,3$.
\endlemma

\proof
If $\sec_2\supset\{l_1,\ldots,l_7\}$
(disjoint from~$m_1$, $m_2$, and each other), then the root
\[*
l_1+\ldots+l_7+3e_2+e_3-h-e_1-2m_1-2m_2
\]
separates $e_1$ and~$e_3$.
\endproof

\lemma\label{lem.trig.sec}
If $\sec_1$ and $\sec_2$ are as in \autoref{cor.trig}\iref{trig.discrete},
then ${\sec_1}*{\sec_2}\le1$.
\endlemma

\proof
If $\sec_1\supset\{l_1,l_2\}$ and $\sec_2\ni m$ so that
$m\cdot l_1=m\cdot l_2=1$, then the root
\[*
e_1+e_2+e_3+l_1+l_2+m-h
\]
separates~$e_1$ and $e_3$ (as well as $l_1$, $l_2$).
\endproof

\lemma \label{lem:triangular.pencil}
For a triangular geometric graph~$\graph$ one has either
\[*
\graphex\same\Psi_{42},\Psi_{38},\Psi'_{36},\Psi''_{36},\ \text{\rm or}\ \Psi^1_{36}
\]
\rom(see \autoref{rem.ext} and \autoref{tab.list}\rom)
or $\ls|\graph|\le35$.
\endlemma

\proof
We proceed as above, using \autoref{cor.trig} and Lemmata~\ref{lem.trig.trig}
and~\ref{lem.trig.sec}. Since we always assume that
\[*
\ls|\sec_1|\ge\ls|\sec_2|\ge\ls|\sec_3|,
\]
we never need to deal with a set $\sec_i$ as in
\autoref{cor.trig}\iref{trig.trig}.

First, we analyze, as in \cite[Appendix~D]{degt.Rams:octics},
all graphs $\graph\supset\pencil$ with $\ls|\pencil|\ge18$, arriving at
$\ls|\graph|\le35$ with the five exceptions listed in the statement.

Then, we use \cite[Appendix~C]{degt.Rams:octics} to list all graphs
$\graph\supset\fiber$ with $\ls|\sec\fiber|\ge19$: we find
that $\ls|\graph|\le35$
unless $\graphex\same\Psi_{42}$, $\Psi_{38}$, or $\Psi''_{36}$.

In the remaining cases,
$\ls|\graph|=\ls|\pencil|+\ls|\sec\fiber|\le17+18=35$.
\endproof


\subsection{Graphs of other types}\label{s.other}
Other types of graphs exhibit no new phenomena: we proceed literally as
in~\cite{degt.Rams:octics} and arrive at essentially the same collection of
graphs as in the smooth case, see~\cite{degt:lines}. Therefore, we merely
state the ultimate result.


\lemma\label{lemma:easygraphs}
Let $\graph$ be a geometric graph of type $\fiber\not\same\tA_2$ or~$\tA_3$.
Then\rom:
\roster*
\item
for $\fiber\same\tA_4$, either
$\graphex \same \Phi_{30}', \Phi_{30}'', \Phi_{30}^{5\prime\prime}$
\rom(see \cite[Lemma 7.7]{degt.Rams:octics}\rom) or
$\ls|\graph|\le29$\rom;
\item
for $\fiber\same\tD_4$, one has $\ls|\graph|\le28$\rom;
\item
for $\fiber\same\tA_5$, either
$\graphex \same \Lambda^\mathrm{A}_{24},\Lambda^4_{24}$
\rom(see \cite{degt:lines} and \autoref{rem.Lambda}\rom) or
$\ls|\graph|\le23$\rom;
\item
for $\fiber\same\tD_5$, either
$\graphex \same \Lambda_{24},\Lambda'_{24},\Lambda''_{24}$
\rom(see \cite{degt:lines}\rom) or $\ls|\graph|\le23$.
\endroster
In all other cases, $\ls|\graph|\le23$.
\endlemma

It is worth mentioning that the somewhat weaker inequality $\ls|\graph|\le29$
for any \emph{locally elliptic} (\ie, $\fiber$-graph with $\mu(\fiber)\ge5$)
geometric graph~$\graph$ and \emph{any} degree~$h^2$ is immediate without a
machine aided computation, see \cite[\S\,7.1]{degt.Rams:octics}.

\remark\label{rem.Lambda}
The hexagonal graph $\Lambda^\mathrm{A}_{24}$ is the one introduced in
\cite[Addendum~4.7]{degt:lines}, whereas
$\Fano_6(\Lambda^4_{24})\supset\Fano_6(\Lambda^\mathrm{A}_{24})$ is an
index~$3$ extension with the same set of lines and four nodes. It is
interesting that, in the case of octics, instead of $\Lambda^4_{24}$ we have an
index~$4$ extension
$\Fano_8(\Lambda^6_{24})\supset\Fano_8(\Lambda^\mathrm{A}_{24})$ with six
nodes, see \cite[\S\,7.1]{degt.Rams:octics}.
\endremark

\subsection{Proof of \autoref{th.main}}\label{proof.main}
Let $\graph$ be a geometric graph and $\ls|\graph|\ge36$. From
\cite[(4.3)]{degt.Rams:octics} we infer that $\graph$ is hyperbolic,
whereupon, arguing type-by-type,
Lemmata~\ref{lem:tA3biquadranglefree}, \ref{lem:tA3type1and2},
\ref{lem.111}, \ref{lem:triangular.pencil}, and~\ref{lemma:easygraphs}
(see also \autoref{rem.ext}) imply that $\graphex$ is one of the exceptions
(other than $\Theta^1_{36}$) listed in \autoref{tab.list}.
By \autoref{th.geometric}, these are all extended graphs of $36$ or more
lines on $K3$-sextics $X\to\Cp4$
\emph{provided that $\NS(X)\otimes\Q$ is generated by
the polarization~$h$ and the classes of lines.}

It remains to analyze, for each of the six graphs~$\graph$ of rank
$\rank\graph<20$, its corank~$1$ extensions by an extra exceptional
divisor~$e$, \cf. \cite[\S\,3.3.2]{degt.Rams:octics}. This is straightforward
for the graphs $\graph\same\Phi'_{36},\Phi''_{36},\Theta_{36},\gf1{36}$ of
cardinality~$36$ and almost straightforward for $\graph\same\Phi_{38}$ (where
we can ``loose'' up to $2$ vertices): by \autoref{prop.exceptional}, the
support $\|e\|\subset\graph$ is an independent set of cardinality at
most~$6$, and we apply the line-by-line algorithm of
\autoref{s.line.by.line}.
The graph $\Psi_{42}$, where up to $6$ lines can be lost, needs extra work;
we deal with it by choosing an appropriate basis for $\NS(X)\otimes\Q$.
Ultimately, the only graph admitting a corank~$1$ extension is $\Theta_{36}$:
we obtain $\Theta_{36}^1$.
\qed

\remark \label{rem:whatthehell}
Comparison of maximal extended Fano graphs
(see \cite{degt:lines,Degtyarev.Rams.triangular,degt.Rams:octics})
shows that the same large configurations of lines keep
appearing on \KK-surfaces of various degrees,
even though the surfaces themselves
clearly differ. It is natural to ask
if the ubiquity of certain configurations
has a geometric meaning.
\endremark

\remark \label{rem:lem.111}
Observe that,
the only pair of extended Fano graphs in \autoref{tab.list} that
share the same plain graph is $\Theta_{36}$ and $\Theta^1_{36}$.
Note that the exceptional divisor in $\Theta^1_{36}$ is disjoint from all
lines, meaning that the incidence graph of lines in~$X$ equals that
of their images in $X_6\iin\Cp4$.
\endremark


\section{Certain automorphisms of smooth sextic $K3$-surfaces}\label{S.aut}

\newcommand{\Pl}{\Pi}
\newcommand{\prj}{\pi}
\newcommand{\coherent}{concordant}

In this section we discuss the notion of deck translation
defined
by a line (\autoref{s.tau}), introduce \emph{special \KK-sextics}
and give their characterization in terms of $3$-isotropic vectors
(see \autoref{s.special}). Then we
characterize
pairs of deck translations generating finite dihedral groups
and
prove
\autoref{cor.aut} in \autoref{proof.cor.aut}.

\subsection{The deck translation defined by a line}\label{s.tau}
Recall that a smooth $K3$-sextic $X\iin\Cp4$ is a
complete
intersection $V_2\cap V_3$ of a quadric and a cubic; the former is
uniquely determined
by~$X$. Fix a line $\ell\subset X$ and consider a plane $\Pl\supset\ell$.
If $\Pl\not\subset V_2$, then $V_2\cap \Pl=\ell+l$ is a degree~$2$ curve and either
\roster
\item\label{i.2.points}
$V_3\cap l$ consists (generically) of three points one of which is in~$\ell$, or
\item\label{i.line}
$V_3\supset l$, in which case $l$ is a line on~$X$ such that $l\cdot\ell=1$.
\endroster
(Since $X$ is smooth, $X\cap \Pl$ has no multiple components.)
If $\Pl\subset V_2$, then
\roster[\lastitem]
\item\label{i.conic}
$V_3\cap \Pl=\ell+c$, where $c$ is a conic on~$X$, possibly reducible;
$c\cdot\ell=2$.
\endroster
In case~\iref{i.conic} we say that $(c,\ell)$ is a \emph{conic-line pair} or that $c$
is the \emph{residual conic} of~$\ell$ (for the uniqueness of the residual conic see \autoref{rem-unique3iso} below).
Thus, as is well known (see \eg, \cite[\S\,5]{Saint-Donat}),
the projection $\Cp4\dashrightarrow\Cp2$ from~$\ell$
restricts to a generically two-to-one
map $\prj_\ell\:X\to\Cp2$. This map is ramified over a sextic curve
$C_\ell\subset\Cp2$, which is smooth apart from a number of nodes, \viz. the
blow-downs of the lines~$l$ as in~\iref{i.line} and conics~$c$ as
in~\iref{i.conic}. The latter may degenerate to cusps, which happens when the
corresponding conic splits into a pair of lines.

The deck translation $\tau_\ell\:X\to X$ of the ramified double
covering~$\prj_\ell$ is an automorphism of~$X$. This automorphism does \emph{not} preserve
the original polarization~$h$: instead, it preserves the
hyperelliptic polarization
$h-\ell=\prj_\ell^*[\Cp1]$.

We are interested in the action $(\tau_\ell)_*$
of~$\tau_\ell$ on the lattices $\NS(X)\subset H_2(X)$. Let
\[*
\rho_a\:H_2(X)\to H_2(X),\quad
x\mapsto x\mp(x\cdot a)a,\quad
a\in H_2(X),\ a^2=\pm2,
\]
be the reflection. Then $(\tau_\ell)_*$ is the composition of
the reflection $-\rho_{h-\ell}$
and the
product of (pairwise commuting) local
contributions of all singular points of the ramification
locus~$C_\ell$. These local contributions are
\roster*
\item
$\rho_e$ for a node~$\bA_1$ with exceptional divisor~$e$, or
\item
$\rho_{e_1+e_2}$ for a cusp~$\bA_2$ with exceptional divisors~$e_1,e_2$.
\endroster
Thus,
\[
(\tau_\ell)_*=-\rho_{h-\ell}\circ\prod\rho_{l_i}\cdot\rho_c,
\label{eq.tau}
\]
where $c$ is the residual conic of~$\ell$, if present, and $l_i$ runs over all lines
in $\Star\ell$ \emph{that are not components of~$c$}, \ie,
the $\bA_1$-summands in \autoref{th.star}.
Note that this expression does not depend on whether the residual conic~$c$, if
present, is irreducible or reducible.

\remark\label{rem.deck.translation}
This
description of the deck translation should be well known to the experts,
but we could not immediately find a convenient reference.
Consider (the minimal resolution of singularities of) a double plane
$\prj\:X\to\Cp2$ ramified over an even degree curve with simple singularities
only.
The projection $\prj\:X\to\Cp2$ is known to
factor through the minimal embedded resolution of singularities
\[*
\prj\:X\overset q\longto Y\longto\Cp2
\]
of the ramification locus, and the
transfer $q^*$ establishes
an isomorphism
\[*
q^*\:H_2(Y;\Q)\overset\same\longto\Ker\bigl[(1-\tau_*)\:
 H_2(X;\Q)\to H_2(X;\Q)\bigr].
\]
Thus, it suffices to ``guess'' an involution that acts correctly on the
exceptional divisors on~$X$ and extends \via\ $-\id$ on their orthogonal
complement. This map is $-\rho_h$, $h:=p^*[\Cp1]$, composed with the
contribution $-\tau_S$ of each singular point~$S$:
\roster*
\item
$\tau_S$ is induced by the only non-trivial symmetry of the Dynkin diagram for
a point $S$ of type $\bA_{>1}$, $\bD_\text{odd}$, or $\bE_6$, and
\item
$\tau_S=\id$ for all other simple singularities.
\endroster
For nodes~$\bA_1$ and cusps~$\bA_2$ this is expressed in terms of
reflections as above.
\endremark


\subsection{Special \vs. non-special sextics}\label{s.special}
A line $\ell\subset V_2$ lies in a plane $\Pi\subset V_2$ if and only if the
quadric $V_2$ is singular. It follows that the existence of a residual conic
of a line $\ell\subset X$
is a property of~$X$ itself rather than that of an individual line.

\definition\label{def.special}
A smooth $K3$-sextic $X\iin\Cp4$ is called \emph{special} if some
(equivalently, each) line $l\subset X$ has a residual conic.
Alternatively, $X=V_2\cap V_3$ is special if and only if
the quadric $V_2$ is singular. In this latter form, we extend the definition
to $K3$-sextics with singularities.
\enddefinition

We have a simple homological characterisation of special sextics in the spirit of
that for octics,
see~\cite{Saint-Donat} or \cite[Lemma~2.27]{degt.Rams:octics}.
To our surprise, we could not find a reference in the literature:
classically, the
$3$-admissiblity in
the context of generators of
the defining ideals
is discussed in
\cite[\S\,7]{Saint-Donat}, where the author explicitely assumes the linear
system
to be \emph{of degree at least eight} --- see
\cite[(7.1)]{Saint-Donat}.

\theorem\label{th.special}
A $K3$-sextic $X\tto\Cp4$ is non-special if and only if the $6$-polarized
lattice $\NS(X)\ni h$ is $3$-admissible
\rom(see \autoref{def.admissible}\rom).
\endtheorem

We precede the proof of \autoref{th.special} with a lemma. Observe that, if $p\in\NS(X)$ is a
$3$-isotropic vector, then so is $h-p$; we call these vectors
\emph{complementary}.

\lemma\label{lem.special}
Let $X$ be a $K3$-sextic and $p,p'\in\NS(X)$ two $3$-isotropic vectors.
\roster
\item\label{spec.2p}
If $X$ is smooth, then either
\roster*
\item
$p\cdot p'=0$, and then $p=p'$, or
\item
$p\cdot p'=3$, and then $p+p'=h$.
\endroster
\item\label{spec.line}
If $l\subset X$ is a line, then $l\cdot p\in\{0,1\}$.
\item\label{spec.conic}
If $c\subset X$ is a conic, then $c\cdot p\in\{0,1,2\}$.
\endroster
\endlemma

\proof
In fact, the statement holds for any $2$-admissible $6$-polarized lattice
$N\ni h$, where the smoothness is understood as
condition~\iref{i.exceptional} in \autoref{def.admissible}.

Let $s:=p\cdot p'$. Then (\cf. \autoref{obs.ker}),
\[*
\det(\Z h+\Z p+\Z p')=6s(3-s)\ge0,\quad\text{hence $0\le s\le3$}.
\]
The values $s=0,3$ result in $\det=0$, implying a relation.
For $s=1$ (resp.\ $s=2$), $e:=p-p'$ (resp.\ $h-p-p'$) is
an exceptional divisor.

Similarly, $\det(\Z h+\Z p+\Z l)\ge0$ implies $-1\le l\cdot p\le2$, and the
value $l\cdot p=-1$ (resp.\ $2$) results in a $2$-isotropic vector $p-l$
(resp.\ $h-p-l$).

Finally, $\det(\Z h+\Z p+\Z c)\ge0$ implies $-1\le c\cdot p\le3$, and the
value $c\cdot p=-1$ (resp.\ $3$) results in a $1$-isotropic vector $p-c$
(resp.\ $h-p-c$).
\endproof

\proof[Proof of \autoref{th.special}]
If $V_2$ is singular, it contains a plane. The latter  meets $V_3$
along a plane cubic curve~$E$. The class $[E]\in\NS(X)$ is
$3$-isotropic.

Conversely, \autoref{lem.special}\iref{spec.line} and \iref{spec.conic}
allow us to  repeat \latin{verbatim}
the proof of \cite[Observation~2.23]{degt.Rams:octics} and show that a $3$-isotropic vector can be chosen in the
closure of
$\Nef X$. Then \cite[(2.26)]{degt.Rams:octics} implies the existence of
an elliptic curve $E\subset X$ such that $h\cdot E = 3$. Finally, the curve $E$  must
be a plane cubic and the plane (projectively) spanned by~$E$ is contained in~$V_2$,
so that the latter cannot be smooth.
%
\endproof

\remark \label{rem-unique3iso}
Assume that $X$ is special and \emph{smooth}, so that there are but two
complementary $3$-isotropic vectors $p_1,p_2\in\NS(X)$.
Then, by \autoref{lem.special}\iref{spec.2p} each line $\ell\subset X$
has a unique orthogonal (\ie, such that $p\cdot\ell=0$) $3$-isotropic vector
$p:=p(\ell)$ and, hence, a \emph{unique} residual conic
$c:=c(\ell)=p(\ell)-\ell$. (This explains the
``the'' in the definition after \autoref{s.tau}\iref{i.conic};
\cf. also \autoref{th.star} and \autoref{rem.conics}.)
This conic, characterized by the identity $c\cdot\ell=2$,
may be reducible, $c=\ell_1+\ell_2$, in which case
$\ell,\ell_1,\ell_2$ constitute a triangle in the Fano graph $\Fn X$.

Conversely, a conic
$c\subset X$, irreducible or reducible,
is part of a (unique) conic-line pair if and only if
$c\cdot p_i\in\{0,2\}$ (as opposed to $c\cdot p_1=c\cdot p_2=1$). The line
$\ell\subset X$ in the pair is $\ell:=p_i-c$, where $p_i\cdot c=0$; it is
characterized by $\ell\cdot c=2$.
\endremark

Let $\ell_1,\ell_2\subset X$ be two lines and $c_i:=c(\ell_i)$, $i=1,2$,
their residual conics. We say that $\ell_1,\ell_2$ are \emph{\coherent}~ if
$p(\ell_1)=p(\ell_2)$; otherwise, $p(\ell_1)+p(\ell_2)=h$ and the lines
$\ell_1,\ell_2$ are called \emph{complementary}.

In the former case
(\coherent), \emph{assuming that $\ell_1$ is not a component of~$c_2$},
since $p(\ell_1)\cdot p(\ell_2)=0$ and all intersections are
non-negative,
we have
\[
\gathered
\ell_1\cdot\ell_2=\ell_1\cdot c_2=c_1\cdot\ell_2=c_1\cdot c_2=0\rlap{\quad and}\\
0\le\ell\cdot(\ell_1+c_1)=\ell\cdot(\ell_2+c_2)\le1
\endgathered
\label{eq.coherent}
\]
for any \emph{other} (\ie, neither $\ell_1,\ell_2$ nor a component of
$c_1,c_2$, if any) line~$\ell$.
In the latter case (complementary),
$\ell_1+c_1+\ell_2+c_2$ is a hyperplane section of~$X$ and
\[
\gathered
\ell_1\cdot(\ell_2+c_2)=\ell_2\cdot(\ell_1+c_1)=1,\quad
c_1\cdot(\ell_2+c_2)=c_2\cdot(\ell_1+c_1)=2,\\
\text{and any other line intersects
(at~$1$)
exactly one of
$\ell_1,c_1,\ell_2,c_2$}.
\endgathered
\label{eq.complementary}
\]

If $\ell_1\cdot\ell_2=0$, the two lines can be
\coherent~or
complementary. If $\ell_1\cdot\ell_2=1$, the lines must be complementary
unless they are two sides of a triangle $\ell_1,\ell_2,\ell$, in which case
all three
lines
are \coherent. The following statement is immediate.

\corollary\label{lem.bipartite}
If a smooth $K3$-sextic~$X$ is special, then the Fano graph $\Fn X$ is either
triangular or bipartite\rom; in the latter case, the
bipartition is given by
 the intersection $\ell\cdot p\in\{0,1\}$ with a fixed $3$-isotropic
vector $p\in\NS(X)$.
\done
\endcorollary

Thus, if  $X$ is
smooth and
special and  $\Fn X$ is triangle free, then $\ls|\mbox{Fn} X| \le 24$
as each part of the bipartition has cardinality at most $12$ by
\autoref{th.Kummer}. This last bound is sharp, see \autoref{rem.Humbert.spec}
below.

Each tiangular Fano graph~$\graph$ is obviously special. In this case,
still assuming~$X$ smooth,
we have a decomposition
$\graph=\pencil_1\cup\pencil_2$, $\pencil_1\cap\pencil_2=\varnothing$.
Each $\pencil_i$ is either discrete, and then $\ls|\pencil_i|\le12$, or a
type~$\tA_2$ pencil, and then $\pencil_{3-i}=\sec\graph_i$ and either
$\pencil_i\same6\tA_2\oplus3\bA_1$ or $\ls|\pencil_i|\le20$, see
\cite[Corollary 3.8]{degt:lines}. In particular, $\graph$ has at most $12$
triangles. For example, $\Psi_{42}=\pencil_1\cup\pencil_2$ with
$\pencil_1\same\pencil_2\same6\tA_2\oplus3\bA_1$.

\proof[Proof of \autoref{cor.boundsforspecial}]
The graphs $\Psi_{42}$, $\Psi_{38}$, $\Psi^1_{36}$ are triangular
(see \autoref{lem:triangular.pencil}),
so they are realized by special \KK-sextics. This implies the bounds for
special \KK-sextics and shows that no non-special \KK-sextic contains more than $36$ lines.
Examples with $36$ lines are given by
\autoref{rem:discussion.graphs}(4).

As just explained a smooth special \KK-sextics contains at most $24$
 lines, and an example of $24$ lines is given in \autoref{rem.Humbert.spec}.
\endproof


\subsection{Proof of \autoref{th.aut}}\label{proof.aut}
We pick a pair of lines $\ell\i,\ell\ii\subset X$ and consider the
correponding deck translations $\tau\i,\tau\ii$.

In view of \eqref{eq.tau},
the actions of $\tau\i_*,\tau\ii_*$ on the lattices $\NS(X)\subset H_2(X)$
depend on the respective stars $\Star\ell\i$, $\Star\ell\ii$ and on
a contribution from the residual
conics (which do not need to split into lines)
when
$X$ is special.
Furthermore, it suffices to restrict the action to the $\Q$-vector
space generated by~$h$, $\ell\i$, $\ell\ii$, and
$\Star\ell\i\cup\Star\ell\ii$: both deck translations act identically on the
orthogonal complement thereof. Thus, our immediate goal is a description
of the \emph{joint star} $(\Star\ell\i\cup\Star\ell\ii)$:  this is
done in \autoref{s.aut.trig}--\autoref{s.aut.skew} below on the case-by-case
basis.

 The following observations are the starting point of each description:
\roster*
\item
the description of an individual star $\Star\ell^*$ given by \autoref{th.star};
\item
properties of complete biquadrangles~\eqref{eq.complete} given by
\autoref{cor.biquad};
\item
properties of residual conics, see mostly~\eqref{eq.coherent}
and~\eqref{eq.complementary}.
\endroster
We begin with a $\{1,2\}$-colored graph~$\graph_0$ of \emph{relevant} lines and
conics, check if it is subgeometric (see \autoref{def.graph}), and,
if necessary (in \autoref{s.aut.skew}), proceed as in \autoref{s.line.by.line}
to extend~$\graph_0$ to a subgeometric graph~$\graph$ by
adding a few more relevant lines. Then we merely determine the order of
the composition $\tau\ii_*\circ\tau\i_*$ on the lattice $\Fano(\graph)$.
It is
crucial that we deal with $\Fano(\graph)=(\Z\graph+\Z h)/\!\ker$ rather than
with $\Z\graph+\Z h$ itself, on which, formally, involutions~\eqref{eq.tau}
also act: in some cases, a non-trivial action on the latter descends to a
trivial one on the former.

We discuss the composition $\tau\ii_*\circ\tau\i_*$ in
\autoref{s.aut.trig}--\autoref{s.aut.skew} below, depending on the original
pair $\ell\i,\ell\ii$;
Lemmata~\ref{lem.aut.intr} and~\ref{lem.aut.skew}
imply \autoref{th.aut}.
\qed

\convention
To save space, in \autoref{s.aut.trig}--\autoref{s.aut.skew} below,
when describing the initial configuration $\graph_0$, we assume
that all
intersections
not mentioned explicitly vanish.
\endconvention

\newcommand{\lB}{m}
\newcommand{\lC}{l}

\subsection{Two sides of a triangle}\label{s.aut.trig}
Assume that the two lines $\ell\i,\ell\ii$ are two sides of a triangle
$\ell\i,\ell\ii,\ell$;
in other words, $X$ is special, $\ell\i\cdot\ell\ii=1$, and $\ell\i,\ell\ii$
are \coherent.
Then, apart from~$\ell$, the joint star consists of
certain lines
\[
\lB_1\i,\ldots,\lB_s\i,\lC_1\i,\ldots,\lC_t\i,\lB_1\ii,
\ldots,\lB_s\ii,\lC_1\ii,\ldots,\lC_r\ii,
\quad
\lB_i\i\cdot \lB_i\ii=1,\ i=1,\ldots,s,
\label{eq.star.trig}
\]
with  $\lB_{i}\i, \lC_{i}\i \in  \Star\ell\i$  (resp.\
 $\lB_{i}\ii, \lC_{i}\ii \in  \Star\ell\ii$).

Let $\fd:=\tau\ii_*\circ\tau\i_*-\id$ and consider the
vector
\[*
v:=12\ell-4\sum_{i=1}^{s}(\lB\i_i+\lB\ii_i)-3\sum_{i=1}^{t}\lC\i_i-3\sum_{i=1}^{r}\lC\ii_i.
\]
By direct computation 
 $\fd^3(v)=0$ and
$\fd^2(v)=\kappa(\ell\i,\ell\ii)(h-\ell\i-\ell\ii-\ell)$, where
\[*
\kappa(\ell\i,\ell\ii):=72-8s-6t-6r.
\]
Starting from the \latin{a priori} bound $0\le r\le t\le9-s\le9$
given by \autoref{th.star}, we find
$104$ subgeometric configurations, all satisfying
$\kappa(\ell\i,\ell\ii)\ge6$.
It follows that
$\tau\ii_*\circ\tau\i_*$ has a
Jordan block (with the eigenvalue $\Gl=1$) of size at least~$3$.
Therefore,
\[ \label{eq:concordantmeet}
\mbox{if \coherent~ lines } \ell\i,\ell\ii
\mbox{ meet, then } \order(\tau\ii_*\circ\tau\i_*)=\infty \, .
\]


\newcommand{\lan}{r}
\newcommand{\lmm}{t}
\newcommand{\lrr}{q}
\newcommand{\laa}{n}
\newcommand{\ldd}{c}  
\subsection{A pair of intersecting lines}\label{s.aut.intr}
Assume that $\ell\i\cdot\ell\ii=1$ but $\ell\i,\ell\ii$ do not intersect a
common third line. Then, in addition to~\eqref{eq.star.trig}, the joint star
may contain a certain number $\lrr\le4$ of common complete biquadrangles
\[
\laa\i_{i1},\laa\i_{i2},\laa\ii_{i1},\laa\ii_{i2},\quad
\laa\i_{ij}\cdot \laa\ii_{ik}=1,\ i=1,\ldots,r,\ j,k=1,2.
\label{eq.star.intr}
\]
Besides, if $X$ is special (and then $\ell\i,\ell\ii$ are complementary),
there also are
\[*
\text{two conics $\ldd\i,\ldd\ii$ (possibly split) such that
$\ldd\i\cdot\ell\i=\ldd\ii\cdot\ell\ii=\ldd\i\cdot \ldd\ii=2$},
\]
(the linear components of $\ldd\i, \ldd\ii$, if any, are not included into
$\lB_*^*,\lC_*^*$).
As in \autoref{s.aut.trig}, we put $\fd:=\tau\ii_*\circ\tau\i_*-\id$ and consider the
vector
\[*
v:=6\ell\ii-2\sum_{i=1}^{s}\lB\i_i-3\sum_{i=1}^{\lmm}\lC\i_i.
\]
One has $\fd^3(v)=0$ and
$\fd^2(v)=\kappa(\ell\i,\ell\ii)(h-\ell\i-\ell\ii)$, where
\[*
\kappa(\ell\i,\ell\ii):=72-12\lrr-8s-6\lmm-6\lan-24\Gd_X
\]
and $\Gd_X=1$ (resp.\ $\Gd_X=0$)
if $X$ is special (resp.\ otherwise).
The \latin{a priori} bounds are $\lrr,s\ge0$ and $0\le \lan\le \lmm\le9-2\lrr-s$; we
find $300$ subgeometric configurations.
We have
\roster*
\item
$\kappa(\ell\i,\ell\ii)\ge0$ for any pair of lines $\ell\i,\ell\ii$ as in
this section;
\item
$\kappa(\ell\i,\ell\ii)=0$ if and only if $X$ is special, $s=0$, and
$(\lrr,\lmm,\lan)$ is one of
\[*
(0, 4, 4), (0, 6, 2), (0, 8, 0);
(1, 4, 2), (1, 6, 0);
(2, 2, 2), (2, 4, 0);
(3, 2, 0);
(4, 0, 0);
\]
\item
in each of these exceptional cases, $\tau\i_*$, $\tau\ii_*$ commute;
hence, so do $\tau\i$, $\tau\ii$.
\endroster
Combining this with the results of \eqref{eq:concordantmeet},
we arrive at the following statement.

\lemma\label{lem.aut.intr}
Assume that two lines $\ell\i,\ell\ii\subset X$ intersect
{\rm(}\ie, $\ell\i\cdot\ell\ii=1${\rm)}. Then
the involutions $\tau\i,\tau\ii$ commute if and only if
\roster
\item\label{aut.intr.special}
the sextic $X$ is special,
\item\label{aut.intr.complementary}
the lines $\ell\i,\ell\ii$ are complementary, and
\item\label{aut.intr.params}
parameters $(\lrr,s,\lmm,\lan)$ in~\eqref{eq.star.trig} and~\eqref{eq.star.intr}
satisfy the restrictions
\[*
s=0,\quad
\lrr,\lmm,\lan\ge0,\quad
2\lrr+\lmm+\lan=8,\quad
\text{\rm$\lmm$ and~$\lan$ are even}.
\]
\endroster
In all other cases, $\tau\i$ and $\tau\ii$ generate an infinite dihedral
group.
\done
\endlemma


\subsection{A pair of skew lines}\label{s.aut.skew}
Apart from $\ell\i,\ell\ii$, there may be certain lines
\[
\lB_1,\ldots,\lB_s,\lC_1\i,\ldots,\lC_{\lmm}\i,\lC_1\ii,\ldots,\lC_{\lan}\ii,
\quad
\lB_i\cdot\ell\i=\lB_i\cdot\ell\ii=\lC_j\i\cdot\ell\i=\lC_k\ii\cdot\ell\ii=1.
\label{eq.star.skew}
\]
We have $0\le s\le3$, see \autoref{th.biquad}, and $0\le \lan\le \lmm\le9-s$, see
\autoref{th.star}.
Besides, if $X$ is special, there are also
\[*
\text{two conics $\ldd\i,\ldd\ii$ (possibly split) such that
$\ldd\i\cdot\ell\i=\ldd\ii\cdot\ell\ii=2$}.
\]
We distinguish between the two cases, see \autoref{s.special}
and~\eqref{eq.coherent}, \eqref{eq.complementary}:
\roster
\item\label{aut.d0}
$\ell\i,\ell\ii$ are \coherent:
$\ldd\i\cdot\ell\ii=\ldd\ii\cdot\ell\i=\ldd\i\cdot \ldd\ii=0$ and
$c_j\i\cdot \ldd\ii=\lC_k\ii\cdot \ldd\i=1$;
\item\label{aut.d1}
$\ell\i,\ell\ii$ are complementary:
$\ldd\i\cdot\ell\ii=\ldd\ii\cdot\ell\i=\ldd\i\cdot \ldd\ii=1$, $s=0$,
$\lC_*^*\cdot \ldd^*=0$.
\endroster
(The linear components of the conics $\ldd\i,\ldd\ii$, if any, are not
listed as lines
$\lB_*^*,\lC_*^*$.)
In case~\iref{aut.d1},
we also need to distinguish between irreducible and
reducible conics. Indeed, if, say, $\ldd\ii$ is irreducible, it does not
contribute to $\tau\i$, whereas if $\ldd\ii=d\ii_1+d\ii_2$ splits so that
$d\ii_1\cdot\ell\i=1$ and $d\ii_2\cdot\ell\i=0$,
then $\tau\i_*$ has an extra reflection defined by~$d\ii_1$.

\figure
\hbox to\hsize{\hss\vbox{%
\setcounter{enumi}0\relax
\csname @hyper@itemtrue\endcsname
\def\itm{\refstepcounter{enumi}%
 \expandafter\xdef\csname @currentlabel\endcsname{\theenumi}\rm(\theenumi)\enspace}
\def\gp#1{,\ \DG{#1}}
\ialign{\hss$#$\hss&&\quad\hss$#$\hss\cr
\tconfig
\. \. \. \. \cr
\. \. \. \. \cr
\. \. \. \. \cr
\. \. \. \. \cr
\endtconfig&
\tconfig
\1 \. \. \. \. \cr
\. \1 \. \. \. \cr
\. \. \1 \. \. \cr
\. \. \. \1 \. \cr
\. \. \. \. \1 \cr
\endtconfig&
\tconfig
\1 \. \. \. \. \. \cr
\. \1 \. \. \. \. \cr
\. \. \1 \. \. \. \cr
\. \. \. \1 \. \. \cr
\. \. \. \. \1 \. \cr
\. \. \. \. \. \1 \cr
\endtconfig&
\tconfig
 \1 \1 \. \. \. \. \. \cr
 \1 \. \1 \. \. \. \. \cr
 \1 \. \. \1 \. \. \. \cr
 \. \. \. \. \1 \. \. \cr
 \. \. \. \. \. \1 \. \cr
 \. \. \. \. \. \. \1 \cr
 \. \. \. \. \1 \1 \1
\endtconfig&
\cr\noalign{\medskip}
\itm\label{DG.4}s\le3\gp4&
\itm\label{DG.5}s=0\gp6&
\itm\label{DG.6}s\le3\gp6&
\itm\label{DG.7.1}s\le2\gp6\cr\noalign{\bigskip}
\tconfig
\. \. \. \. \. \. \. \cr
\. \1 \1 \. \. \. \. \cr
\. \1 \. \1 \. \. \. \cr
\. \. \1 \1 \. \. \. \cr
\. \. \. \. \1 \1 \. \cr
\. \. \. \. \1 \. \1 \cr
\. \. \. \. \. \1 \1
\endtconfig&
\tconfig
 \1 \1 \1 \. \. \. \. \. \cr
 \1 \. \. \1 \1 \. \. \. \cr
 \1 \. \. \. \. \. \. \. \cr
 \. \1 \. \. \. \. \1 \. \cr
 \. \. \1 \. \. \. \1 \. \cr
 \. \. \. \1 \. \. \. \1 \cr
 \. \. \. \. \1 \. \. \1 \cr
 \. \. \. \. \. \1 \1 \1
\endtconfig&
\tconfig
\1 \1 \. \. \. \. \. \. \cr
\1 \. \1 \. \. \. \. \. \cr
\. \1 \. \1 \. \. \. \. \cr
\. \. \1 \1 \. \. \. \. \cr
\. \. \. \. \1 \1 \. \. \cr
\. \. \. \. \1 \. \1 \. \cr
\. \. \. \. \. \1 \. \1 \cr
\. \. \. \. \. \. \1 \1
\endtconfig&
\tconfig
\1 \1 \1 \. \. \. \. \. \. \cr
\1 \. \. \1 \1 \. \. \. \. \cr
\1 \. \. \. \. \1 \. \. \. \cr
\. \1 \. \1 \. \. \1 \. \. \cr
\. \1 \. \. \. \. \. \1 \. \cr
\. \. \1 \. \. \1 \. \1 \. \cr
\. \. \. \1 \. \. \. \. \1 \cr
\. \. \. \. \1 \1 \. \. \1 \cr
\. \. \. \. \. \. \1 \1 \1
\endtconfig\cr\noalign{\medskip}
\itm\label{DG.7.2}s\le2\gp6&
\itm\label{DG.8=6}s\le1\gp6&
\itm\label{DG.8=8}s=0\gp8&
\itm\label{DG.9}s=0\gp6\cr
\crcr}}\hss}
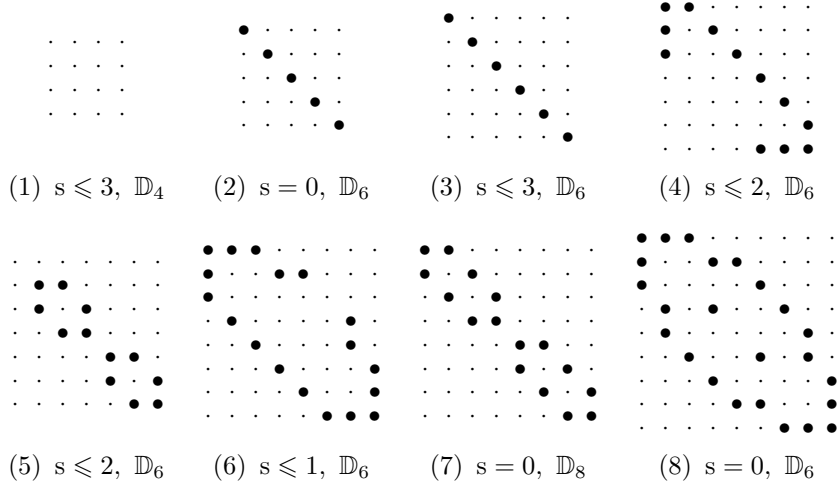
\caption{Pairs generating finite dihedral groups (see \autoref{lem.aut.skew})}\label{fig.finite}
\endfigure

The problem
is that we cannot describe in simple terms the intersections
\smash{$\lC\i_i\cdot \lC\ii_j$}: \latin{a priori} we can only assert that each $\lC\i_i$
intersects at most three \smash{$\lC\ii_j$}'s and \latin{vice versa},
see \autoref{th.biquad}. Therefore, we
start from a collection of lines $\lB_i$, \smash{$\lC\i_j$} and conics $\ldd\i,\ldd\ii$ (if
present, possibly also split into lines) and add the remaining lines
\smash{$\lC\ii_j$}
one by one, as in \autoref{s.line.by.line},
obtaining $9749$ configurations. The ``geometric insight'' for
the algorithm described in
\autoref{s.line.by.line} is that each support $\|\lC_j\ii\|$,
see~\eqref{eq.supp},
has size at
most~$3$ and $\fm=0$ in~\eqref{eq.fm}.

Unlike \autoref{s.aut.trig} and \autoref{s.aut.intr},
we could not find a simple ``universal'' explanation of the fact that
$\tau\ii_*\circ\tau\i_*$ has infinite order: sometimes the matrix has a
non-trivial Jordan block, but more often its characteristic polynomial has an
irreducible factor that fails to be cyclotomic. There are but $22$ exceptions,
listed in the next statement.


\lemma\label{lem.aut.skew}
Assume that two lines $\ell\i,\ell\ii\subset X$ are skew,
$\ell\i\cdot\ell\ii=0$. Then
the involutions $\tau\i,\tau\ii$ generate a finite dihedral group if and
only if $\Star\ell\i\cup\Star\ell\ii$ is as shown in \autoref{fig.finite} and
\roster*
\item
in case~\iref{DG.4}, $X$ is special and $\ell\i,\ell\ii$ are \coherent\rom;
here, $\tau\i,\tau\ii$ commute\rom;
\item
in case~\iref{DG.5}, $X$ is special, $\ell\i,\ell\ii$ are complementary,
and $\ldd\i,\ldd\ii$ do \emph{not} split\rom;
\item
in case~\iref{DG.6}, either $X$ is non-special \rom($0\le s\le3$\rom) or
$X$ is special and $\ell\i,\ell\ii$ are
\coherent~\rom(only for $s=0,1,3$\rom)\rom;
\item
in all other cases, $X$ is necessarily non-special.
\endroster
Shown in the figure are
the intersection pattern of lines \smash{$\lC_j\i,\lC_k\ii$},
the possible number $s$ of common lines~$\lB_i$ in~\eqref{eq.star.skew},
and the group $G$ generated by $\tau\i,\tau\ii$.
\done
\endlemma

\remark\label{rem.aut.skew}
Note that
the full configuration of lines
may be
larger than
shown in \autoref{fig.finite}. For example,
in \autoref{fig.finite}\iref{DG.6}
with~$X$
special, all twelve lines $\lC\i_j,\lC\ii_k$ are \coherent; hence, there \emph{must} be
at least six
extra lines $\lC_p$, $p=1,\ldots,6$, such that
\smash{$\lC_p\cdot \lC\i_p=\lC_p\cdot \lC\ii_p=\lC\i_p\cdot \lC\ii_p=1$}
(and $\lC_p\cdot \lC_*^*=0$
otherwise); in particular, the full graph $\Fn X$ is necessarily triangular.
That is why we speak about \emph{subgeometric} graphs: we only take into
account lines and conics relevant for $\tau\i,\tau\ii$.
\endremark

\subsection{Proof of \autoref{cor.aut}}\label{proof.cor.aut}
Assume that,
for each pair $\ell\i,\ell\ii$ of lines on~$X$, the deck
translations $\tau\i,\tau\ii$ generate a finite group.
Step-by-step we conclude that
\roster*
\item
$\Fn X$ contains no triangles by \eqref{eq:concordantmeet},
\item
there is a pair of adjacent vertices in $\Fn X$ by \autoref{lem.aut.skew},
\item
$X$ is special by \autoref{lem.aut.intr}\iref{aut.intr.special}, and
\item
any pair $\ell\i,\ell\ii\subset X$ of \coherent~skew lines is as in
case~\iref{DG.4} of \autoref{fig.finite}
by \autoref{lem.aut.skew} and
\autoref{rem.aut.skew}.
\endroster
We assert that $\Fn X$
has no complete biquadrangles; then, due to
\autoref{lem.aut.intr}\iref{aut.intr.params}, it also has no quadrangles.
Indeed, given~\eqref{eq.complete}, each pair $l_i,l_j$,
$1\le i<j\le3$, must be as in \autoref{fig.finite}\iref{DG.4} and, hence,
there must be twelve pairwise disjoint (\eg, because they are
 \coherent,
see~\eqref{eq:concordantmeet})
lines $s_{ij}$,
$1\le i\le3$, $1\le j\le4$, such that $s_{ij}\cdot l_i=1$ and,
by \autoref{cor.biquad}, $s_{ij}$ is
disjoint from all other lines in~\eqref{eq.complete}. Together with
$m_1,m_2,m_3$ this makes $15$ pairwise disjoint lines, contradicting
\autoref{th.Kummer}.

Now, by Lemmata~\ref{lem.aut.intr}\iref{aut.intr.params} and \ref{lem.aut.skew},
there is a line $\ell_0\in\Fn X$ adjacent to at least four, necessarily
pairwise disjoint, lines $\ell_1,\ldots,\ell_4$. Each pair $\ell_i,\ell_j$,
$1\le i<j\le4$, is as in \autoref{fig.finite}\iref{DG.4}, giving rise to
eight
extra lines $s_{i1},\ldots,s_{i4},s_{j1},\ldots,s_{j4}$ as in the figure.
All $16$ lines $s_{ij}$, $1\le i,j\le4$, are \coherent\
and, by~\eqref{eq:concordantmeet}, pairwise disjoint,
contradicting \autoref{th.Kummer}.
\qed

\section{Humbert sextics}\label{S.Humbert}

The \emph{Humbert cubic line complex}~\cite{Humbert} is the line
complex~$\Cal{C}$ cut
by a certain cubic hypersurface on the quadric
$\operatorname{Gr}(4,2)\subset\Cp5$.
Studied in~\cite{DDK} are the so-called \emph{Humbert sextic
$K3$-surfaces}, which are the transversal hyperplane sections of~$\Cal{C}$; this
is indeed a $5$-parameter family of smooth $K3$-sextics.


\begin{table}
\begin{center}
\caption{The incidence relation between $12+12$ lines \smash{(see \cite{DDK})}}\label{tab.Humbert}
\hrule height0pt
\def\rb#1{\rotatebox[origin=c]{270}{\ $#1$\ }}
\def\cb#1{\hbox to11pt{\hss#1\hss}}
\def\vr{\vrule width0pt height11pt}
\def\*{\rlap{$^*$}}
\def\*{}
\begin{tabular}{|c|cccc|cccc|cccc|}\hline
&\vr\cb{1\*}&\cb{2\*}&\cb{3}&\cb{4}&\cb{5}&\cb{6\*}&\cb{7}&\cb{8\*}&\cb{9\*}&\cb{10}&\cb{11}&\cb{12\*}\\
 \hline\vr
1\*&$\bullet$&$\bullet$&&&&$\bullet$&&$\bullet$&$\bullet$&& &$\bullet$ \\
2&$\bullet$&$\bullet$&&&$\bullet$&&$\bullet$&&&$\bullet$&$\bullet$&\\
3&&&$\bullet$&$\bullet$&&$\bullet$& &$\bullet$&&$\bullet$ &$\bullet$&\\
4&&&$\bullet$&$\bullet$&$\bullet$&&$\bullet$&&$\bullet$&&&$\bullet$\\
 \hline\vr
5&&$\bullet$&&$\bullet$&$\bullet$&$\bullet$&&&  &$\bullet$&&$\bullet$\\
6&$\bullet$&&$\bullet$&&$\bullet$&$\bullet$&&&$\bullet$&&$\bullet$&\\
7&&$\bullet$&&$\bullet$&&&$\bullet$&$\bullet$&$\bullet$&&$\bullet$&\\
8&$\bullet$&&$\bullet$&&&&$\bullet$&$\bullet$&&$\bullet$&&$\bullet$\\
 \hline\vr
9&$\bullet$&&&$\bullet$&&$\bullet$ &$\bullet$&&$\bullet$&$\bullet$&&\\
10&&$\bullet$&$\bullet$&&$\bullet$&&&$\bullet$&$\bullet$&$\bullet$& &\\
11&&$\bullet$&$\bullet$&&&$\bullet$ &$\bullet$&&&&$\bullet$&$\bullet$\\
12&$\bullet$&&&$\bullet$&$\bullet$&&&$\bullet$&&&$\bullet$&$\bullet$ \\
\hline\omit
\end{tabular}%
\end{center}
\end{table}
A very general Humbert sextic~$X$ has two $12$-tuples
of pairwise skew lines, \viz.
\[
\text{\emph{$\Ga$-lines} $l_1,\ldots,l_{12}$\quad and\quad
 \emph{$\Gb$-lines} $m_1,\ldots,m_{12}$}
\label{eq.Humbert.ab}
\]
cut on~$X$ by some $\Ga$- and $\Gb$-planes in the Grassmannian;
the
intersections $l_i\cdot m_j$ are as shown in \autoref{tab.Humbert}.
The set $\Ga:=\{l_1,\ldots,l_{12}\}$ of $\Ga$-lines splits into three
\emph{quartets}
\[*
\Ga_1:=\{l_1,\ldots,l_4\},\quad\Ga_2:=\{l_5,\ldots,l_8\},\quad
 \Ga_3:=\{l_9,\ldots,l_{12}\}
\]
shown in \autoref{tab.Humbert}; likewise, we have
$\Gb:=\{m_1,\ldots,m_{12}\}=\Gb_1\cup\Gb_2\cup\Gb_3$.

We reserve the notation $\Humbert:=\Fn X=\Ga\cup\Gb$
for the Fano graph of a very general Humbert sextic $X\iin\Cp4$.
The graph~$\Humbert$ has $16$ complete biquadrangles, \eg,
\[*
\HL[ 1, 5, 11, 14, 18, 24 ];
\]
each has three $\Ga$-lines, one from each $\Ga_r$, and three $\Gb$-lines, one
from each~$\Gb_s$. Any choice of a pair of lines from two distinct $\Ga$-quartets,
say, $l_i\in\Ga_1$ and $l_j\in\Ga_2$, gives rise to a unique
biquadrangle. The same holds for two $\Gb$-quartets.

\observation[see~\cite{DDK}]\label{rem.quads}
In addition to the $16$ complete biquadrangles,
$\Humbert$ has $18$ \emph{proper quadrangles} (\ie, those not
contained in a complete biquadrangle), two in each $(4\times4)$-\emph{cell}
in \autoref{tab.Humbert} (or, equivalently, pair $(\Ga_r,\Gb_s)$ of quartets,
one from each family);
these quadrangles are clearly seen in the table.
\endobservation

\subsection{Degenerations}\label{s.degenerations}
The goal of this section is the equilinear stratification of the
space of smooth Humbert sextics, see \autoref{th.Humbert}, which was
announced in \cite[Remark 7.1]{DDK}.
We start with a description of  smooth
rational curves of low degree.

\observation[see~\cite{DDK}]\label{rem.curves}
Apart from the $24$ lines, a very general Humbert sextic $X$ has
\roster
\item\label{i.conics}
$9$ conics $c_{rs}$, indexed by the $(4\times4)$-cells
in \autoref{tab.Humbert}, see \autoref{rem.quads};
\item\label{i.cubics}
no twisted cubics;
\item\label{i.quartics}
$72$ quartics $q_{ij}:=h-l_i-m_j$, indexed by the \emph{disjoint} pairs $(l_i,m_j)$.
\endroster
For the conics, one has $c_{rs}=h-\sum q$, where $q$ is any of the two
proper quadrangles contained in the corresponding cell. In other words,
$c_{rs}+q$ is a hyperplane section of~$X$, and $c_{rs}$ is the common conic of two
such sections. One has
\[
c_{rs}\cdot c_{uv}=
 \begin{cases}
  0,&\text{if $u=r$ or $v=s$},\\
  2,&\text{if $u\ne r$ and $v\ne s$}.
 \end{cases}
\label{eq.conics}
\]

Needless to say that these statements, especially the abscence
of other smooth rational curves of low degree, are proved using Vinberg's
\autoref{alg.Vinberg}.
\endobservation

\observation[see~\cite{DDK}]\label{rem.aut.G}
The automorphism group
\[*
\Aut\Humbert\same(\SG4\times\SG4)\rtimes\Z/2
\]
acts
transitively on the set of $16$ complete biquadrangle, so that the stabilizer of a
biquadrangle $q$ is
$\Aut q\same(\SG3\times\SG3)\rtimes\Z/2$.
Alternatively, $\Aut\Humbert$ induces
the full group $(\SG3\times\SG3)\rtimes\Z/2$ of symmetries of
the $(3\times3)$-grid in \autoref{tab.Humbert}
(and it is in this sense that the splitting into quartets is natural),
and the stabilizer of a
$(4\times4)$-cell
$c\same\tA_3\oplus\tA_3$ induces the index~$2$ subgroup of
$\Aut c\same(\DG8\times\DG8)\rtimes\Z/2$ that does not mix $\Ga$-
and $\Gb$-lines.
\endobservation

\theorem\label{th.Humbert}
Any line on a Humbert sextic~$Y$ is either \rom(the limit
of\rom) one of the $24$ original lines~\eqref{eq.Humbert.ab}
or a component of one of the $9$ conics in
\autoref{rem.curves}\iref{i.conics}.
There are eight proper equilinear strata, see
\autoref{fig.degenerations} and \autoref{conv.Humbert}.
\figure
\[*
\hbox to\hsize{\hss\vbox{\ialign{&\hss$#$\hss\cr
\Gr=16\,&\multispan2\hss$\Gr=17$\,\hss&\multispan3\hss$\Gr=18$\,\hss
 &\multispan2\hss$\Gr=19$\,\hss\cr\noalign{\smallskip}
\vrule\
\tconfig
 \1 \. \. \cr
 \. \. \. \cr
 \. \. \. \cr
\endtconfig
\ \vrule\,&
\vrule\
\tconfig
 \1 \1 \. \cr
 \. \. \. \cr
 \. \. \. \cr
\endtconfig
\ \vrule&
\ \tconfig
 \1 \. \. \cr
 \. \1 \. \cr
 \. \. \. \cr
\endtconfig
\ \vrule\,&
\vrule\ \tconfig
 \1 \1 \. \cr
 \1 \. \. \cr
 \. \. \. \cr
\endtconfig
\ \vrule&
\ \tconfig
 \1 \. \. \cr
 \. \1 \. \cr
 \. \. \1 \cr
\endtconfig
\ \vrule&
\ \tconfig
 \1 \1 \. \cr
 \. \. \1 \cr
 \. \. \1 \cr
\endtconfig
\ \vrule\,&
\vrule\ \tconfig
 \1 \1 \. \cr
 \1 \1 \. \cr
 \. \. \. \cr
\endtconfig
\ \vrule&
\ \tconfig
 \1 \1 \. \cr
 \1 \. \1 \cr
 \. \1 \1 \cr
\endtconfig
\ \vrule\cr\noalign{\medskip}
\{1\}\,&\{1\}&\Z/2\,&\Z/2&\SG3&\,(\Z/2)^2&(\Z/2)^2&\,\DG{12}\cr
\crcr}}\hss}
\]
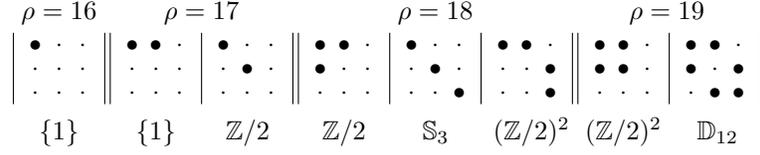
\caption{Smooth degenerations of Humbert sextics (see \autoref{th.Humbert})}\label{fig.degenerations}
\endfigure
\endtheorem

In the last stratum, the one with six split conics, all sextics have $36$ lines;
their Fano graphs are isomorphic to $\Theta_{36}$ in \autoref{tab.list}.

\convention\label{conv.Humbert}
In \autoref{fig.degenerations}, we indicate
\roster*
\item
the Picard rank $\rho:=\rank\NS(Y)$ and
\item
group $\Aut_h(Y)$ of projective automorphisms
\endroster
of a very general representative~$Y$ of each proper stratum.
Besides, we show the conics $c_{rs}$
\rom(equivalently, cells of the $(3\times3)$-grid in
\autoref{tab.Humbert}\rom) that split into pairs of new lines.
Part of the statement of \autoref{th.Humbert} is that this set of split conics,
regarded up to the action of $\Aut\Humbert$, see
\autoref{rem.aut.G}, determines the stratum.
\endconvention

The following consequence was used in the proof of \cite[Lemma~5.1]{DDK}.

\corollary\label{cor.Humbert.aut}
The graph $\Fn Y$ of any Humbert sextic~$Y$ has exactly two maximal
independent subsets of size~$12$, \viz. $\Ga$ and $\Gb$. Hence, both
$\Aut(\Fn Y)$ and $\Aut_hY$ preserve the original $24$ lines $\Ga\cup\Gb$.
\done
\endcorollary


\proof[Proof of \autoref{th.Humbert}]
As explained in \autoref{s.line.by.line},
an extra line $u\in\Fn Y$ is determined by its
support~$\|u\|$, see~\eqref{eq.supp}, and,
by \autoref{cor.biquad}, this set $\|u\|$ must have exactly one common vertex
with each
complete biquadrangle $q\subset\Humbert$. Using the transitivity of $\Aut\Humbert$,
we can assume that $l_1\in\|u\|$; then
\[*
\|u\|\sminus l_1\subset
 \Humbert\sminus Q_1=\HL[ 2, 3, 4, 15, 16, 17, 19, 22, 23 ],\quad
 Q_1:=\bigcup q,
\]
the union running over the four complete biquadrangles containing~$l_1$.
Furthermore, $\|u\|$ must have \emph{exactly one} common vertex with the
difference $q\sminus Q_1$ for each of the other $12$ complete biquadrangles
$q\not\ni l_1$. Assuming that
$l_2\in\|u\|$ (again,
without loss of generality) and
arguing as above, we find that $\|u\|$ must have exactly one common vertex
with
each of the four sets
\[*
\HL[ 3, 15 ],\ \HL[ 3, 16 ],\
  \HL[ 4, 15 ],\ \HL[ 4, 16 ].
\]
Thus, up to $\Aut\Humbert$, the support is
\[*
\|u\|=\HL[1,2,3,4]\quad\text{or}\quad\|u\|=\HL[1,2,15,16].
\]
The former case $\|u\|=\Ga_1$ 
contradicts
\autoref{th.Kummer}: the union $\Gb\cup u$ would
be a collection of $13$ pairwise disjoint lines.

\remark\label{rem.4h}
According to \cite[(2.11)]{DDK}, one has $\sum\Humbert=4h$ in $\NS(X)$.
Thus, it is not an accident that the support of an extra lines consists of
exactly four
vertices.
\endremark

In the latter case, $\|u\|=\HL[1,2,15,16]$, the graph $\graph':=\Humbert\cup u$
is subgeometric and its saturation contains another line $v:=c_{11}-u$ (\cf.
\autoref{rem.curves});
it has
$\|v\|=\HL[3,4,13,14]$.
We have $u\cdot c_{11}=v\cdot c_{11}=-1$ and $u\cdot v=1$;
geometrically, this means that the conic $c_{11}$ is no longer irreducible:
it splits into $u+v$.

In view of \autoref{th.biquad}, in any geometric configuration
$\graph'\supset\Humbert$ there may be at most one line intersecting both~$l_1$
and~$l_2$ (and the same holds for $l_3,l_4$); hence, a pair $u,v$ of extra
lines is indeed uniquely determined by the conic~$c_{rs}$ that splits or,
equivalently, pair $(\Ga_r,\Gb_s)$ of quartets such that
$\|u\|,\|v\|\subset\Ga_r\cup\Gb_s$; these new lines, denoted by
$\ell'_{rs},\ell''_{rs}$,
extend
each of the two proper
quadrangles in the $(r,s)$-th cell of \autoref{tab.Humbert} (see
\autoref{rem.quads})
to complete biquadrangles.

It remains to analyze the combinations of conics that can split
simultaneously. In addition to a collection $\Cal{P}\subset\{1,2,3\}^2$ of
pairs of indices (a set of supports),
we need to specify the matrix~$\fm$ in~\eqref{eq.fm}.
In view of~\eqref{eq.conics}, we have
\[*
\ell'_{rs}\cdot\ell'_{uv}=\ell'_{rs}\cdot\ell''_{uv}
 =\ell''_{rs}\cdot\ell'_{uv}=\ell''_{rs}\cdot\ell''_{uv}=0
\]
whenever $u=r$ or $v=s$; otherwise,
\[*
\ell'_{rs}\cdot\ell'_{uv}=\ell''_{rs}\cdot\ell''_{uv}=1,\qquad
 \ell'_{rs}\cdot\ell''_{uv}=\ell''_{rs}\cdot\ell'_{uv}=0
\]
\emph{under an appropriate ordering of the new lines}.
Remarkably, arguing as in \autoref{s.line.by.line}, we find that, for
each $\Cal{P}$, there is a unique,
up to $\Aut\Humbert$, choice of $\fm$
that makes the new graph $2$-admissible. Furthermore,
for
\[*
\Cal{P}=\tconfig
 \1 \1 \1 \cr
 \. \. \. \cr
 \. \. \. \cr
\endtconfig\quad\text{or}\quad
\Cal{P}=\tconfig
 \1 \1 \. \cr
 \1 \1 \. \cr
 \. \. \1 \cr
\endtconfig
\]
the resulting configuration is not geometric. It remains to select the
configurations that are saturated (see \autoref{def.graph})
and establish the irreducibility of each stratum
as explained in \autoref{s.moduli}.
\endproof

\remark\label{rem.Humbert.spec}
A very general Humbert sextic~$X$ is non-special; however, there is a
$4$-parameter family of special Humbert sextics. For proof, one merely
considers the lattice
\[*
N:=(\Z\Humbert+\Z h+\Z p_\Ga)/\!\ker,
\]
where
\[*
p_\Ga^2=0,\quad p_\Ga\cdot h=3,\quad p_\Ga\cdot l_i=1,\quad p_\Ga\cdot m_i=0
\]
for $i=1,\ldots,12$; the last two intersections are dictated by the
bipartition of~$\Humbert$, see \autoref{lem.bipartite}.
Then, it is easy to
see that $N$ is geometric and $\Fn(N,h)=\Humbert$.
In a sense, we extend our machinery to weighted graphs and let
$v^2=\operatorname{wt}(v)\in\Z$.

None of the proper strata in \autoref{fig.degenerations} (see
\autoref{th.Humbert})
contains
a special sextic: the Fano graphs are neither triangular nor bipartite,
see \autoref{lem.bipartite}.

Apart from Humbert, there are ten other equilinear families
of smooth sextics with a
bipartite graph of $12+12$ lines. Only one of them has a codimension~$1$
stratum (depending on $2$ parameters) of special sextics. Proofs are
based on the line-by-line algorithm of \autoref{s.line.by.line}; they
will appear
elsewhere.
\endremark

\remark\label{rem.Humbert.curves}
A very general representative~$Y$ of each of the proper equilinear strata in
\autoref{fig.degenerations}
(see \autoref{th.Humbert})
has a great deal of twisted cubics and extra quartics (\cf.
\autoref{rem.curves}).
With two exceptions, $Y$ also has a few extra conics. Proofs based on
Vinberg's \autoref{alg.Vinberg} are left to the reader.
\endremark

\subsection{Digression: the equiconical strata}\label{s.equiconic}
Similar to~\autoref{s.degenerations}, one can analyze the strata of Humbert
sextics with
an extra conic: one needs to consider a bi-colored graph $\Humbert\cup u$ with an
extra vertex~$u$ of degree $u\cdot h=2$.

If $u\cdot v=2$ for some $v\in\Humbert$, then $u$ is the residual conic of
the line~$v$ and, hence, the sextic is special, see \autoref{rem.Humbert.spec}.
Therefore, we can assume that $u\cdot v\in\{0,1\}$ for each line $v\in\Humbert$
and argue as in~\autoref{s.line.by.line}:
the support $\|u\|\subset\Humbert$ is an $8$-element set (see
\autoref{rem.4h}) with exactly two common points with each complete
biquadrangle (\cf. \autoref{cor.biquad}). Up to $\Aut\Humbert$ this leaves but
the five possibilities below.

\medskip\noindent
\underline{Case 1}: $\|u\|=\HL[ 1, 2, 3, 4, 13, 14, 15, 16 ]$.
This orbit consists of
the original nine conics that are present in any Humbert sextic,
see \autoref{rem.conics}\iref{i.conics}.

\medskip\noindent
\underline{Case 2}: $\|u\|=\HL[ 1, 2, 3, 4, 5, 6, 19, 20 ]$.
This orbit of length $72$ results in graphs that are not extensible: the root
\[*
\def\HLcomma{+}
\def\HLform#1{#1}
-h+u+\HL[ 1, 3, 18, 20 ]
\]
separates~$l_5$ and~$m_5$.

\medskip\noindent
\underline{Case 3}: $\|u\|=\HL[ 1, 2, 3, 4, 5, 6, 7, 8 ]=\Ga\sminus\Ga_3$.
This $6$-element orbit consists of all sets of the form $\Ga\sminus\Ga_r$ and
$\Gb\sminus\Gb_s$, $r,s=1,2,3$. It results in a $4$-parameter family of
Humbert sextics
containing all six extra conics.

\medskip\noindent
\underline{Case 4}: $\|u\|=\HL[ 1, 2, 3, 5, 12, 15, 19, 21 ]$.
The orbit consists of $96$ elements and results in a $4$-parameter family of
Humbert sextics. A very general member has six extra conics, \eg, $u$ itself,
\[*
\|v\|=\HL[ 4, 5, 9, 10, 11, 13, 20, 23 ],\quad
\|w\|=\HL[ 4, 6, 7, 8, 12, 14, 18, 22 ],
\]
and the three supports obtained by the automorphism $l_i\leftrightarrow m_i$
of $\Humbert$.

\medskip\noindent
\underline{Case 5}: $\|u\|=\HL[ 1, 2, 5, 6, 15, 16, 19, 20 ]$.
The orbit consists of $72$ elements and gives rise to a $4$-parameter family
of Humbert sextics; a very general member has $12$ extra conics, all from the
same orbit. This family is studied in detail in
\cite[\S\S\,5.3 and~5.4]{DDK}; it is the only known $4$-parameter family of
Humbert sextics admitting a non-trivial (anti-symplectic) projective automorphism.

\medskip
Thus, the last three cases give us three $4$-parameter families of
non-special Humbert
sextics with extra conics; in each case, a very general member has no extra
lines.
Unlike the case of extra lines, each family contains a $3$-parameter
subfamily of special sextics.
One could continue and analyze further degenerations with more extra
conics and/or lines but, since the present paper is mainly about lines, we do
not pursue this line of research.

\subsection{The ramification loci}\label{s.sextics}
Fix a very general Humbert sextic~$X$ and a line $\ell\subset X$, \eg,
$\ell=l_1$, consider the projection $\pi_\ell\:X\to\Cp2$, see
\autoref{s.tau}, and let $C\subset\Cp2$ be its ramification locus.
The sextic curve~$C$ is described in \cite[\S\,3]{DDK}.
The eight lines
\[*
l_5,l_{11};\ l_6,l_{9};\ l_7,l_{10};\ l_8,l_{12}
\]
project, in
pairs, to the four sides of a complete quadrilateral $Q\subset\Cp2$, and $C$
has nodes at the six vertices of~$Q$: they are the images of the six $\Gb$-lines
\[*
m_1,m_2,m_6, m_8, m_9,m_{12}
\]
that intersect~$\ell$. The lines $l_2,l_3,l_4$ are mapped to the diagonals
of~$Q$; the images are tangent to~$C$ at some smooth points.
Finally, the remaining
$\Gb$-lines
\[*
m_3,m_4,m_5,m_7,m_{10},m_{11},
\]
\viz. those disjoint from~$\ell$, are mapped to generic tritangents to~$C$.
In other words, there is a complete quadrilateral $Q\subset\Cp2$ such that
\[
\text{$C$ has nodes at the vertices of~$Q$ and is tangent to the diagonals of~$Q$};
\label{eq.C-quad}
\]
besides, $C$ has six extra tritangents.
We are interested in the extent to which these restrictions
determine~$C$. The following statement is announced in \cite[\S\,6]{DDK}, where some
parts are proved geometrically. Given~$C$, we denote by
$Y\to\Cp2$ the minimal resolution of singularities of the double plane
ramified over~$C$. It is a $K3$-surface.

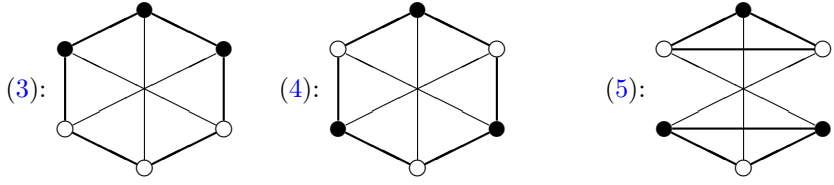
\begin{figure}
\unitlength1.5pt
\def\w{\circle{4}}
\def\ww{\circle{6}}
\def\b{\circle*{4}}
\def\tl{\thicklines}
\def\hhex{%
\put(2,9){\tl\line(2,-1){16}}
\put(38,9){\tl\line(-2,-1){16}}
\put(2,31){\tl\line(2,1){16}}
\put(38,31){\tl\line(-2,1){16}}
\put(20,2){\line(0,1){36}}
\put(2,11){\line(2,1){36}}
\put(38,11){\line(-2,1){36}}
}
\def\hex{%
\hhex
\put(0,12){\tl\line(0,1){16}}
\put(40,12){\tl\line(0,1){16}}
}
\def\lb#1{\put(-15,18){\eqref{#1}:}}
\begin{picture}(40,40)(0,0)
\hex
\put(20,0)\w
\put(0,10)\w
\put(40,10)\w
\put(20,40)\b
\put(0,30)\b
\put(40,30)\b
\lb{fBl.6}
\end{picture}
\qquad\qquad
\begin{picture}(40,40)(0,0)
\hex
\put(20,0)\w
\put(0,10)\b
\put(40,10)\b
\put(20,40)\b
\put(0,30)\w
\put(40,30)\w
\lb{fBl.2}
\end{picture}
\qquad\qquad\qquad
\begin{picture}(40,40)(0,0)
\hhex
\put(2,10){\tl\line(2,0){36}}
\put(2,30){\tl\line(2,0){36}}
\put(20,0)\w
\put(0,10)\b
\put(40,10)\b
\put(20,40)\b
\put(0,30)\w
\put(40,30)\w
\lb{fBl.4}
\end{picture}
\caption{The pull-backs of the diagonals of~$Q$ (see \autoref{th.ramification})}\label{fig.hex}
\end{figure}

\theorem\label{th.ramification}
There are two irreducible $6$-parameter families of sextic curves
$C\subset\Cp2$ satisfying~\eqref{eq.C-quad}. The families differ topologically, by
the adjacency of the pull-backs in~$Y$ of the diagonals of~$Q$\rom:
the six lines in~$Y$ constitute either
\roster
\item\label{fB.good}
a hexagon, as in \autoref{fig.hex}, left \rom(see also \autoref{conv.ramification}\rom), or
\item\label{fB.bad}
two triangles, as in \autoref{fig.hex}, right.
\endroster
Furthermore, there are three irreducible $5$-parameter families of pairs
$(C,\ell)$, where $C$ is a sextic satisfying~\eqref{eq.C-quad} and
$\ell\subset\Cp2$ is a tritangent to~$C$\rom: generically, either
\roster[\lastitem]
\item\label{fBl.6}
  $C$ is as in case~\iref{fB.good} and has six tritangents, or
\item\label{fBl.2}
  $C$ is as in case~\eqref{fB.good} and has two tritangents, or
\item\label{fBl.4}
  $C$ is as in case~\eqref{fB.bad} and has four tritangents.
\endroster
\endtheorem

As explained in \cite[\S\,6]{DDK}, it is case~\iref{fBl.6} that corresponds
to the ramification locus of a
planar model of a Humbert sextic.

\convention\label{conv.ramification}
The pull-back of each diagonal $d_i$ of~$Q$ splits into a pair of lines
$d_i',d_i''\subset Y$ so that $d_i'\cdot d_i''=1$. In \autoref{fig.hex},
pairs $d_i',d_i''$ are depicted by pairs of antipodal vertices connected by
the diagonals of the hexagon. The other six edges
(thick) of the graph constitute
either a hexagon (left) or a pair of triangles (right).

The black vertices are the lines adjacent to a component of the pull-back of
one of the
tritangents, see the proof of \autoref{th.ramification}.
\endconvention

\def\gr#1{\graph_{\ref{#1}}}
\def\ns#1{N_{\ref{#1}}}

\proof[Proof of \autoref{th.ramification}]
Since the double plane~$Y$ is hyper-elliptic, we are interested in
(sub-)geometric $1$-admissible lattices and extensible graphs.

In accordance with our
general paradigm, we start from \emph{lines} only, in the hope that the
exceptional divisors will appear as a result of the computation.
Thus, we consider the graph~$\graph$ consisting of
\roster*
\item
eight pairwise disconnected vertices $l_i',l_i''$, $i=1,\ldots,4$, to play the
r\^{o}le of the pull-backs of the sides of the quadrilateral~$Q$, and
\item
six vertices $d_k',d_k''$, $k=1,2,3$ (pull-backs of the diagonals)
disconnected from all $l_i^*$ and connected pairwise: $d_k'\cdot d_k''=1$.
\endroster
Referring to geometry,
we add six extra edges so that $d_1',d_2',d_3',d_1'',d_2'',d_3''$ is either
\roster
\item\label{d.6}
a single $6$-cycle, in the order listed (see \autoref{fig.hex}, left), or
\item\label{d.3.3}
a pair of $3$-cycles $(d_k')$ and $(d_k'')$ (see \autoref{fig.hex}, right).
\endroster
This gives us two graphs $\gr{d.6}$, $\gr{d.3.3}$ and two
polarized lattices
$N_i:=\Fano_2(\graph_i)$,
and we
analyze their finite index extensions.

During the computation,
we use the shorthand $\bd_i:=d_1'+d_i''$ and $\bl_i:=l_i'+l_i''$.
In both cases, we have $\rank N_i=14$: the radical
$\ker(\Z\graph_i+\Z h)$ is generated by
\[*
2h-2(\bd_1+\bd_2+\bd_3)+\bl_1+\bl_2+\bl_3+\bl_4.
\]

\remark
At first sight, the graphs look too symmetric to reveal the geometric
structure, even the action of the deck translation. This is indeed so (\cf.
the parasitic $2$-nodal family for $\gr{d.3.3}$ below): the symmetry
is broken by
the (unavoidable) passage to a finite index extension.
We use \autoref{alg.Nikulin}.
\endremark

In case~\eqref{d.6}, letting $\graph:=\gr{d.6}$ and $N:=\ns{d.6}$, we have
\[*
\discr N\same(\Z/2)^{10}\oplus(\Z/5),
\]
so that $N$ is not geometric. There are four
$(\Aut\graph)$-orbits of isotropic elements, \viz. those represented by
\[
u:=\frac12(\bd_1-\bd_3),\quad
v:=\frac12(\bl_1-\bl_2),\quad
w:=\frac12(d_1'+d_3'-\bl_1),\quad
u+v;
\label{eq.orbits.1}
\]
the first three contain separating roots, \viz. the representatives shown.
Hence,
the only, up to $\Aut\graph$,
geometric finite index extension of $N$ is that by
$e_1,e_4$,
where
\[
2e_1=\bd_1+\bd_2-\bl_1-\bl_2,\qquad
2e_4=\bd_1+\bd_3-\bl_1-\bl_3;
\label{eq.ee}
\]
by \autoref{lem.extensible} and \autoref{alg.Vinberg},
it has the $14$ original lines and six nodes
\[
\alignedat3
&e_1,\qquad &e_2&:=h-e_1-\bd_3,\qquad &e_3&:=h-e_1-e_4-\bl_1,\\
&e_4,\qquad &e_5&:=h-e_4-\bd_2,\qquad &e_6&:=e_1+e_4-\bd_1+\bl_1,
\endalignedat
\label{eq.nodes}
\]
corresponding to case~\iref{fB.good} in \autoref{th.ramification}.

In case~\eqref{d.3.3}, letting $\graph:=\gr{d.3.3}$ and $N:=\ns{d.3.3}$, we have
\[*
\discr N\same(\Z/2)^6\oplus(\Z/8)^2,
\]
and the discriminant form is even, so that $N$ is not geometric.
There are three $(\Aut\graph)$-orbits of
order~$2$ isotropic elements, \viz. $u$, $v$, and $u+v$ in~\eqref{eq.orbits.1}.
The two former contain separating roots, and this fact rules out
order~$4$ isotropic elements as well: each has a multiple in the
orbit of~$u$. Thus, we find two geometric extensions of $N$: the
index~$2$ extension by $e_1$ and index~$4$ extension by
$e_1,e_4$, see~\eqref{eq.ee}.
The latter has $14$ lines and six nodes~\eqref{eq.nodes},
corresponding to case~\iref{fB.bad} of \autoref{th.ramification}. The former
has $26$ lines but only two nodes, $e_1$ and~$e_2$;
hence, it is not related to our original geometric problem.

For the second part of the statement, we keep the geometric kernel generated
by $e_1,e_4$ and extend $\graph_i$ by an extra vertex~$m'$ to play the r\^{o}le
of a component of the pull-back of~$\ell$. Up to $\Aut\graph_i$, we
can assume that $m'$ intersects all $l_i'$ and is disjoint from all $l_i''$.
Besides, $m'$ must intersect exactly one line from each pair $(d_i',d_i'')$;
the geometrically meaningful adjacencies, two in case~\iref{fB.good} and one
in case~\iref{fB.bad}, are shown in black in \autoref{fig.hex}. We find that
all three lattices are geometric and contain six nodes and $26$, $18$, or
$22$ lines. To reconfirm the fact that the extra $2k=12$, $4$, or~$8$ lines
obtained
are the pull-backs of $k$ tritangents to~$C$, we find that they break into $k$ pairs
$(m_i',m_i'')$ such that $m_i'\cdot m_i''=3$.

Finally, we argue as in \autoref{s.moduli} to establish that each
family is irreducible.
\endproof

{
\let\.\DOTaccent
\def\cprime{$'$}
\bibliographystyle{amsplain}

\providecommand{\bysame}{\leavevmode\hbox to3em{\hrulefill}\thinspace}
\providecommand{\MR}{\relax\ifhmode\unskip\space\fi MR }
\providecommand{\MRhref}[2]{%
	\href{http://www.ams.org/mathscinet-getitem?mr=#1}{#2}
}
\providecommand{\href}[2]{#2}

}

\end{document}